\documentclass[conference]{IEEEtran}

\IEEEoverridecommandlockouts
\usepackage{hyperref}
\usepackage{amsmath,amssymb,amsfonts}
\usepackage{algorithmic}
\usepackage{graphicx}
\usepackage{textcomp}
\usepackage{xcolor}
\usepackage{amsthm}
\usepackage{bm}
\usepackage{dsfont}
\usepackage{subfigure}
\usepackage{subcaption}
\usepackage[algo2e,ruled,noend]{algorithm2e}

\newcommand{\M}{\bm M}
\newcommand{\K}{\mathcal K}

\newcommand{\A}{\mathcal A}

\newcommand{\X}{\mathcal X}
\newcommand{\U}{\mathcal U}
\newcommand{\W}{\mathcal W}

\newcommand{\D}{\mathcal D}
\newcommand{\R}{\mathbb R}
\newcommand{\ball}{\bar{\mathcal{B}} }
\newcommand{\norm}[1]{\left\lVert#1\right\rVert}

\newcommand{\one}{\mathds{1}}

\newtheorem{theorem}{Theorem}

\newtheorem{lemma}{Lemma}
\newtheorem{proposition}{Proposition}
\newtheorem{corollary}{Corollary}

\theoremstyle{definition}
\newtheorem{definition}{Definition}
\newtheorem{assumption}{Assumption}

\def\BibTeX{{\rm B\kern-.05em{\sc i\kern-.025em b}\kern-.08em
    T\kern-.1667em\lower.7ex\hbox{E}\kern-.125emX}}
\begin{document}

\title{Online Nonstochastic Control with Convex Safety Constraints}
\author{Nanfei Jiang, Spencer Hutchinson, Mahnoosh Alizadeh}

\maketitle

\begin{abstract}
This paper considers the online nonstochastic control problem of a linear time-invariant system under convex state and input constraints that need to be satisfied at all times. We propose an algorithm called Online Gradient Descent with Buffer Zone for Convex Constraints (OGD-BZC), designed to handle scenarios where the system operates within general convex safety constraints. We demonstrate that OGD-BZC, with appropriate parameter selection, satisfies all the safety constraints under bounded adversarial disturbances. Additionally, to evaluate the performance of OGD-BZC, we define the regret with respect to the best safe linear policy in hindsight. We prove that OGD-BZC achieves $ \tilde{\mathcal{O}}(\sqrt{T})$ regret given proper parameter choices. Our numerical results highlight the efficacy and robustness of the proposed algorithm.
\end{abstract}

\section{Introduction}

Online nonstochastic control has emerged as a powerful paradigm for managing systems in unpredictable and adversarial settings, attracting significant attention due to its robustness and adaptability \cite{agarwal_online_2019, hazan_introduction_2023}. This framework ensures that control solutions can effectively respond to unforeseen disturbances while optimizing a prescribed performance metric.

Recently, there has been growing interest in integrating safety constraints into the online nonstochastic control framework (e.g., \cite{li_online_2021, liu_online_2023, nonhoff_online_2024, zhou_safe_2023}), increasing its relevance for real-world scenarios. For example, autonomous vehicles are required to navigate through varying traffic patterns without collisions, and robotic systems in manufacturing must avoid causing harm to human workers. A key challenge in these applications is that the safety constraints are often nonlinear, creating a need for algorithms compatible with general convex constraints.

In this context, this paper introduces the Online Gradient Descent with Buffer Zones for Convex Constraints (OGD-BZC) algorithm. Building on the work of \cite{li_online_2021}, which handles affine safety constraints, OGD-BZC generalizes the safety constraints to general convex safety sets, ensuring step-wise safety and achieving $\tilde{\mathcal{O}}(\sqrt{T})$ regret under adversarial cost and disturbance. 
This regret bound represents a nontrivial extension from the work of \cite{li_online_2021}. Specifically, the regret analysis in \cite{li_online_2021} heavily relies on the affine structure of the constraints, resulting in a bound that scales with the number of linear constraints and making it challenging to generalize to the convex case. In contrast, we develop a novel analytical approach that (i) provides a sublinear regret bound for problems with general convex constraints and (ii) eliminates the dependence of the regret bound on the number of constraints, thus improving on the existing literature.
\\

\noindent \textbf{Related Work:}
In the following, we highlight related work that studies online control with constraints, general online nonstochastic control, online convex optimization, model predictive control, and safe reinforcement learning.

\textit{Online Control with Constraints:}
There are several other works  that address online nonstochastic control problem with constraints, \cite{li_online_2021, liu_online_2023, nonhoff_online_2024, zhou_safe_2023, zhou_safe_2023-1}. The work by \cite{zhou_safe_2023, zhou_safe_2023-1} considers time-varying system and safety constraints under unpredictable noise. They propose a gradient-based algorithm called Safe-OGD to achieve bounded dynamic regret against any safe linear policy. This paper assumes the existence of a safe controller for each step, regardless of the state at each step, which might not be practical in some problem settings. The authors in \cite{liu_online_2023} conduct an extensive study on online nonstochastic control, focusing on both soft and hard constraints, which are assessed by cumulative constraint violations. The key difference lies in how violations are measured: for hard constraints, violations are always non-negative for each step, whereas for soft constraints, they can be negative if the state and input stay strictly inside the safety constraints. Their algorithm, named COCA (Constrained Online Nonstochastic Control Algorithm), achieves a regret of $\mathcal{O}(T^{2/3})$ and a constraint violation of $\mathcal{O}(T^{-1/3})$ when dealing with hard constraints. While the focus of \cite{liu_online_2023}  is not on achieving zero constraint violation, the authors mention that using a projection-based method similar to \cite{li_online_2021}, their algorithm can also achieve zero anytime violation. This guarantee of zero violation will only be assured with high probability as their proof is contingent on a random parameter ($Q_t$ representing a virtual queue) being below a threshold, which holds only with probability $1 - 1/T$. In comparison, our safety guarantee is deterministic. The work by \cite{nonhoff_online_2024} also considered convex constraints and uses an approach inspired by the robust model predictive control to ensure constraint satisfaction. In their approach, the algorithm is compared with the optimal steady-state as a regret benchmark, which is different from the definition we adopt. Additionally, their cost function is assumed to be strongly convex and with Lipschitz continuous gradients, which is a stronger condition than we use.


\textit{General Online Nonstochastic Control:}
The foundational work introducing the online nonstochastic control framework by \cite{agarwal_online_2019} shows that disturbance-action controllers can achieve $\mathcal{O}(\sqrt{T}\log(T))$ regret with respect to the best linear policy in hindsight. This work is extended by \cite{agarwal_logarithmic_2019, foster_logarithmic_2020} to provide logarithmic regret under the assumption of strongly convex cost functions. Follow up work on online nonstochastic control (with no constraints) include \cite{gradu_adaptive_2023, minasyan_online_2021, simchowitz_making_2020, simchowitz_improper_2020, zhao_non-stationary_2023}.

 \textit{Online Convex Optimization (OCO):}
The goal of OCO is to minimize the cumulative loss over time, when the loss functions are revealed only after decisions are made \cite{hazan_introduction_2023}. Some work also consider OCO with constraints by allowing soft or hard cumulative constraint violation \cite{cao_virtual-queue-based_2018,  mahdavi_trading_2012, neely_online_2017, yu_online_2017, yuan2018online}. The OCO framework is a crucial tool in our algorithm to provide bounded regret guarantees. However, it is important to note that OCO does not consider the system dynamics or any disturbances.

 \textit{Model Predictive Control:} A traditional method for ensuring constraint satisfaction is robust Model Predictive Control (MPC) \cite{mayne_robust_2005, rawlings_model_2017, schwenzer_review_2021}. It typically relies on a predefined model or cost function to predict future states and optimize control inputs. An adversarial setting complicates this because the cost function can change adversarially, making it difficult to plan optimally over the horizon.

 \textit{Safe Reinforcement Learning:}
 Safe reinforcement learning (RL) focuses on learning the optimal policy by interacting with the environment with safety constraints. Many studies have been conducted in this area (e.g.
\cite{fisac_general_2018, garcia_comprehensive_2015, koller_learning-based_2018}). However, most theoretical studies in the Safe RL literature typically assume environments that do not change over time to provide rigorous guarantees. In contrast, our framework adapts to time-varying cost functions through the analysis of policy regret.
\\

\noindent \textbf{Paper Overview:}
This paper is structured as follows: Section \ref{section: problem form} outlines the problem formulation, including essential definitions and assumptions critical for our analysis. In Section \ref{section: preli}, we present preliminaries on online nonstochastic control. Section \ref{section: online alg} details the development of our proposed Online Gradient Descent with Buffer Zone for Convex Constraints (OGD-BZC) algorithm. Section \ref{section: thm} demonstrates the main theoretical contributions of our work. In Section \ref{section: numerical}, we provide several numerical experiments to validate the effectiveness and robustness of the OGD-BZC algorithm in practical scenarios. Finally, in Section \ref{section: conclusion}, we draw conclusions and outline our future directions.

\subsection{Basic Notations}

In this paper, $\norm{\cdot}_1, \norm{\cdot}_2, \norm{\cdot}_{\infty}$ denote the $L_1, L_2, L_{\infty}$ norm for vectors or matrices. We denote by $\one_S$ the indicater function which equals to $1$ if event $S$ is satisfied and $0$ elsewhere. Occasionally, we use $[H]$ for some integer $H$ to abbreviate the index set $\{ 1,2,...,H \}$. 
We use $\ball_{\norm{\cdot}}^n(\Delta)$ to denote the closed ball centered at zero with radius $\Delta > 0$ in an $n$-dimensional Euclidean space. Often, we omit the dimension or norm notation when there is no confusion in the context. To improve clarity, we express bounds using $\mathcal{O}(\cdot)$ to omit constants that do not depend on the horizon $T$. Similarly, $\tilde{\mathcal{O}}(\cdot)$ refers to the same concept while ignoring the $\log$ factors.

\section{Problem Formulation} 
\label{section: problem form}
We consider a Linear Time-Invariant (LTI) system described by the dynamics:
\begin{equation}\label{dynamics}
    x_{t+1} = Ax_t + Bu_t + w_t, \quad \forall t \ge 0,
\end{equation}
where $x_t \in \mathbb{R}^n$ is the state of the system, $u_t \in \mathbb{R}^m$ is the control input, and $ w_t \in \W := \{w  \in \mathbb{R}^n: \norm{w}_{\infty} \le \overline{w} \}$ represents bounded disturbances. $A \in \mathbb{R}^{n\times n}, B \in \mathbb{R}^{n \times m }$ are known system matrices and $\overline{w}$ is known.

Based on the system dynamics, we consider an online nonstochastic control problem with convex safety
constraints. The problem is stated as follows:

\textbf{Online Nonstochastic Control.} At each step $t \in \{0,1,...T\}$, an agent observes the current state $x_t$ and assigns a control input $u_t$. After $u_t$ is implemented, the agent suffers a stage cost $c_t(x_t,u_t)$. The cost function $c_t(\cdot,\cdot)$ is generated adversarially (i.e. arbitrarily) and revealed to the agent after the control input $u_t$ is taken. The system then evolves to the next step by following the system dynamics in \eqref{dynamics}.

\textbf{Safety Constraints.} To encode safety concerns, we require a general convex constraint on both state $x_t$ and control input $u_t$:
\begin{equation}\label{eq: state safety set}
    x_t \in \X, \quad u_t \in \U, \quad \forall t \ge 0,
\end{equation}
where $\X$, $ \U$ are closed convex sets. An algorithm is said to be \textit{safe} if the generated state $x_t \in \X$ and the control input $u_t \in \U$ for all $0\le t \le T$ and for any sequences of disturbances $\{w_t \in \W \}_{t=0}^T$. We set $x_0 = 0$ and require that both $\X$, $\U$ should include the origin. The case when $x_0 \neq 0$ is left for future work. The following are some assumptions and definitions.



\textbf{Performance Metric.} 
As standard in the online nonstochastic control literature, we consider regret with respect to the class of safe and strongly stable linear policies $\mbox{$u_t = -K x_t$}$, for $K \in \mathbb{R}^{m\times n}$, as the comparison benchmark for our online algorithm.  Strong stability is a quantitative notion of stability that encodes the rate of convergence of any stable system as it implies that the system's response to disturbances or initial conditions diminish at a certain guaranteed rate over time (see \cite{cohen_online_2018}). We provide a formal definition next.
 
\begin{definition}[Strong Stability]
    We define a linear controller $u_t = -K x_t$ to be $(\kappa, \gamma)$-strongly stable for $\kappa >1$ and $\gamma \in (0,1]$ if there exist matrices $L$ and $H$ such that $A-BK = HLH^{-1}$, with $\norm{L}_2 \le 1 - \gamma$ and $\max \left \{\norm{K}_2, \norm{H}_2, \norm{H^{-1}}_2 \right \} \le \kappa$. Additionally, we define $\kappa_B := \max\{\norm{B}_2,1 \}$\footnote{In fact, every stable matrix is $\kappa,\gamma$-strongly stable, for some $\kappa, \gamma$ \cite{agarwal_online_2019}.  
}.
\end{definition}

The regret of an online algorithm $\A$, over $T$ steps, is defined as:
\begin{equation}\label{regret definition}
    \text{Reg}_T(\A) = \max_{\{w_k \in \W \}} \left ( \sum_{t=0}^T c_t(x_t^{\A},u_t^{\A}) - \min_{K \in \K} \sum_{t=0}^T c_t(x_t^K,u_t^K) \right ),
\end{equation}
where $x_t^{\A}$, $ u_t^{\A}$ is the state and control input generated by algorithm $\A$, and $x_t^K, u_t^K$ generated by  the best linear policy $K$ from the class of policies $\K$, defined as,
$$
\begin{aligned}
\K := & \left \{ K \in \mathbb{R}^{n\times m} : K \text{ is } (\kappa,\gamma)\text{-strongly stable}, \right. \\
& \left.  \quad x_t^{K} \in \X, u_t^{K} \in \U \quad \forall w_t \in \W,  0 \le t \le T  
\right \}.
\end{aligned}
$$

For an online algorithm $\A$, a performance guarantee is represented by a sublinear bound on the regret.

Given the adversarial nature of the control problem, we next generalize the concept of a “buffer zone” from \cite{li_online_2021} to the convex constraint case, incorporating it as a margin within the safety constraints. This allows the online algorithm sufficient reaction time to ensure safety even in the face of unforeseen disturbances. The size of this margin will be judiciously chosen to ensure favorable regret. To define this margin, we first define two set operations.

\begin{definition}[Shrinkage and Expansion \cite{hutchinson_impact_2023}]
    The $\Delta$-\textit{shrinkage} of a set $\D \subset \R^d$ under norm $\norm{\cdot}$, denoted by $\D_{\Delta}^{\norm{\cdot}}$ is defined as\footnote{ $\D_{\Delta}^{\norm{\cdot}}$ can be also defined using Minkowski subtraction \cite{schneider_convex_2014}. The Minkowski subtraction of sets $A, B \subseteq \R^d$ is defined as $A \ominus B := \{x: x+B \subseteq A \} $. Thus, we can write $\D_{\Delta}^{\norm{\cdot}} = \D \ominus \ball_{\norm{\cdot}} (\Delta)$ for $\Delta > 0$. }
    \begin{equation}
        \D_{\Delta}^{\norm{\cdot}} := \left\{ x \in \D: x + y \in \D, \forall y \in \ball_{\norm{\cdot}} (\Delta)     \right \},
    \end{equation} 
where $\ball_{\norm{\cdot}} (\Delta)$ it the closed ball with radius $\Delta$, under the same dimension and norm. The $\Delta$-\textit{expansion} of a set $\D \subset \R^d$ under norm $\norm{\cdot}$, denoted by $\D_{-\Delta}^{\norm{\cdot}}$ is defined as
    \begin{equation}
        \D_{-\Delta}^{\norm{\cdot}} := 
        \left \{ x+y: x \in \D, y \in \ball_{\norm{\cdot}} (\Delta)  \right \}.
     \end{equation}
\end{definition}

 It is noteworthy that $\D_0^{\norm{\cdot}}$ is trivially defined to be $\D$ itself. For  simplicity of notation, we will omit the norm notation in the superscript throughout this text. Unless otherwise specified, the term $\D_{\Delta}$ should be interpreted as the $\Delta$-shrinkage/expansion of $\D$ under infinity norm (i.e. $\D_{\Delta}^{\norm{\cdot}_{\infty}}$). 

\begin{definition}[Strictly and Loosely Safe Policies]
A policy/algorithm $\A$ is called $\epsilon$-strictly safe for some $\epsilon>0$ if and only if for all $ 0 \le t \le T$, we have $x_t^{\A} \in \X_{\epsilon} $ and $u_t^{\A} \in \U_{\epsilon} $ under any disturbance sequence $\{ w_k \in \W \}_{k=0}^T$, where $\X_{\epsilon}$ and $\U_{\epsilon}$ denote the $\epsilon$-shrinkage of $\X$ and $\U$. 

Similarly, $\A$ is called $\epsilon$-loose safe if and only if for all $ 0 \le t\le T$, $x_t^{\A} \in \X_{-\epsilon} $ and $u_t^{\A} \in \U_{-\epsilon}$ under any disturbances, where $\X_{-\epsilon}$ and $\U_{-\epsilon}$ are the $\epsilon$-expansion of $\X$ and $\U$.

This definition generalizes the idea of ``buffer zone" defined in \cite{li_online_2021} for affine constraints to general convex constraints.  
As we will see, this buffer zone can be helpful to account for nonstochastic disturbances and approximation errors caused by the problem conversion. Next, we state two assumptions.

\end{definition}

\begin{assumption}[Loss Function] \label{assumption: loss function}
We assume that $c_t(x_t, u_t): \mathbb{R}^n \times \mathbb{R}^m \rightarrow \mathbb{R}$ is convex and differentiable with respect to $x_t$ and $u_t$. Furthermore, for any $\norm{x}_2 \le D$ and $\norm{u}_2 \le D$, there exist a constant $G$ such that   $|c_t(x, u)|$, $\norm{\nabla_x c_t(x, u)}_2$, and $\norm{\nabla_u c_t(x, u)}_2$ are bounded by $GD$. 

\end{assumption}

The following assumption specifies the existence of a strictly safe linear policy over an infinite horizon. This assumption is commonly encountered in the constrained optimization and control literature \cite{li_online_2021, boyd_convex_2004, limon_robust_2010}. It ensures that our performance metric, as defined in \eqref{regret definition}, is well-defined. Moreover, the existence of such a policy also guarantees the feasibility of a safe disturbance-action policy, which is the core approach to be used in this paper.

\begin{assumption}[Existence of Strictly Safe Linear Policy] \label{assumption: existency}
There exists $K_{ss} \in \K $ such that the linear feedback controller $u_t = -K_{ss} x_t$ is $\epsilon_*$-strictly safe in infinite horizon for some $\epsilon_* > 0$.   
\end{assumption}

Here, being infinitely-horizon strictly safe means that $x_t^{K_{ss}} \in \X_{\epsilon_*}$ and $ u_t^{K_{ss}} \in \U_{\epsilon_*}$ for all $t \ge 0$ under any disturbance sequence $\{ w_k \in \W \}_{k=0}^{\infty}$.\footnote{We should mention that our Assumption \ref{assumption: existency} differs slightly from the assumption made in \cite{li_online_2021}, where they only require the existence of an $\epsilon_*$-strictly safe linear policy for a \emph{finite} time horizon $T$. However, this could indirectly allow $\epsilon_*$ to depend on $T$ since $\epsilon_*$ could decrease as $T$ increases. As this could impact the result of the regret analysis, we make a slightly stronger assumption of an infinite horizon ${\epsilon_*}$-strictly safe policy. }

\section{Preliminaries} \label{section: preli}
The online policy we will adopt in this work will be within the class of Disturbance-Action Controllers \cite{agarwal_online_2019}, which we formally define next.
\begin{definition}[Disturbance-Action Controller (DAC)]
    Fix a $(\kappa, \gamma)$-strongly stable matrix $K$. A \textit{disturbance-action controller} $\pi(K,\{\M_t\})$ of memory size $H$ is defined as,
    \begin{equation} \label{def: DAC}
        u_t = -K x_t + \sum_{i=1}^H \M_t^{[i]} w_{t-i}, 
    \end{equation}
    where $w_t = 0$ for $t < 0$ and\footnote{The choice of coefficient $a$ in the definition of $\mathcal{M}$ is not very crucial. For the technical simplicity and without loss of generality, we will pick $a = 2\kappa^3$.}
    \begin{equation}
    \begin{aligned}  \label{eq: M matrix}
    \M_t \in \mathcal{M} := & \left \{ (\M_t^{[1]}, ..., \M_t^{[H]}) : \norm{\M_t^{[i]}}_2 \le a(1-\gamma)^{i-1}, \right. \\
     & \quad \left. \M_t^{[i]}  \in \R^{m\times n}, a > 0, \quad \forall i \in [H]   \right \}.
    \end{aligned}
    \end{equation}
\end{definition}
    
    The DAC encodes the behavior of the past $H$-step disturbances. Modeling with the class of DAC has many advantages. For instance, in \cite{hazan_introduction_2023}, it is shown that the DAC class always gives a stabilizing control signal whenever $K$ stablizes the system. 
    In addition, for brevity of notation, we define $\pi(K,\{\M\})$ to be the DAC class with fixed weight, i.e. $\M_t = \M, \forall t \ge 0$.

    The following proposition, which is directly adopted from \cite{agarwal_online_2019}, allows us to divide the state and control for the DAC into two components, one representing the effect of the predicted state of the system  $H$ steps prior to the current step $t$, and the other, referred to as the  surrogate state and controller, capturing the additional state influences  that are not captured by the predictive model and are due to recent disturbances only. This result will be crucial in describing the set of safe DAC policies.

\begin{proposition}[Lemma 4.3 in \cite{agarwal_online_2019}] \label{prop: DAC}
    When applying a DAC $\pi(K,\{\M_t\})$ as shown in \eqref{def: DAC}, the state and control input can be written as: 
    \begin{equation}\label{eq: DAC}
        x_t = A_K^H x_{t-H} + \tilde{x}_t, \quad u_t = -KA_K^Hx_{t-H} + \tilde{u}_t,
    \end{equation} where $A_K = A - BK$ and the superscript $H$ denotes raising the matrix to the $H$-th power, $\tilde{x}_t$ and $\tilde{u}_t$ are referred to as the  surrogate state and input and 
are defined as:
    \begin{equation}
    \label{eq: tilde x u}
    \begin{aligned}
        \tilde{x}_t  &= \sum_{k=1}^{2H} \Psi_k^x (\M_{t-H:t-1})w_{t-k}, \\
        \tilde{u}_t  &= -K\tilde{x}_t + \sum_{i=1}^{H} \M_t^{[i]} w_{t-i}, \\
                     &= \sum_{k=1}^{2H} \Psi_k^u (\M_{t-H:t})w_{t-k},
    \end{aligned}
    \end{equation}
    where $\M_{t-H:t} := \left \{ \M_{t-H},...,\M_t   \right \}$ and the ``disturbance-to-state/response'' matrices $\Psi_k^x (\M_{t-H:t-1}), \Psi_k^u (\M_{t-H:t})$ are defined as:
    $$
    \begin{aligned}
    \Psi_k^x (\M_{t-H:t-1}) = & A_K^{k-1}\one_{(k\le H)} \\& + \sum_{i=1}^{H} A_K^{i-1} B \M_{t-i}^{[k-i]} \one_{(1 \leq k-i \le H)}, \\
    \Psi_k^u (\M_{t-H:t}) = &\M_t^{[k]} \one_{(k \le H)} - K \Psi_k^x(\M_{t-H:t-1}).
    \end{aligned}
    $$

    Specifically, when we apply the same DAC weight matrix from step $t-H$ to $t$, the disturbance-to-state/response matrices will have simpler expressions, denoted as, 
    \begin{equation}\label{psiring}
        \mathring{\Psi}_k^x (\M) = \Psi_k^x(\M,..., \M),\quad \mathring{\Psi}_k^u (\M) = \Psi_k^u(\M,..., \M).
    \end{equation}

    Often, when there is no confusion in the context, we will omit the weight matrix and simply write $\Psi_k^x, \Psi_k^u$ and $\mathring{\Psi}_k^x, \mathring{\Psi}_k^u$.
\end{proposition}

With the definitions we have established, we are now prepared to develop our online algorithm. 
    

\section{The OGD-BZC Algorithm} \label{section: online alg}

This section will introduce our online algorithm through detailed analysis. As an overview, OGD-BZC is a projected online gradient descent algorithm that finds an appropriate weight matrix $\M_t$ at each time step. For a weight matrix $\M_t$ to be appropriate, it needs to help us achieve low regret while always guaranteeing the safety constraints. To attain low regret, OGD-BZC performs gradient descent on an \emph{approximate cost function}, aiming to reduce the actual cost function at each step. For guaranteeing safety constraints, the algorithm projects $\M_t$ onto a \emph{safe policy set}. This ensures that the chosen policy maintains the system's state and control input within safety constraints under any sequence of disturbances. As such, the main two design components of OGD-BZC are 1) construction  of the safe policy set, and 2) construction of the approximate cost function, as discussed next.

\textbf{Constructing the Safe Policy Set:}

Our current safety set is defined based on the state and control input \eqref{eq: state safety set}. To construct a respective \emph{set of safe policies}, we will first utilize Proposition~\ref{prop: DAC} to make a connection between state space and the policy space.

The first step is to decompose the state and control input into two parts as in \eqref{eq: DAC}. If the state is bounded, the first term $A_K^H x_{t-H}$ will decay exponentially with respect to $H$. This observation means that we can always choose \(H\) large enough to force the first term to become smaller than some threshold \(\epsilon_1\), i.e., \(\norm{A_K^H x_{t-H}}_{\infty} \le \epsilon_1\). Consequently, we can require the surrogate state \(\tilde{x}_t\) to stay within the \(\epsilon_1\)-shrinkage of the safety constraint, i.e., \(\tilde{x}_t \in \X_{\epsilon_1}\), ensuring that the true state will satisfy
$$x_t = A_K^H x_{t-H} + \tilde{x}_t \in \ball^n(\epsilon_1) + \X_{\epsilon_1} = (\X_{\epsilon_1})_{-\epsilon_1} \subseteq \X, $$
where the last set inclusion comes from a property of shrinkage and expansion. The details can be found in the Appendix \ref{Appendix: support}, Proposition \ref{proposition: Minkwski subtraction}. The same analysis can also be applied to the control input $u_t$.
 
Thus, we will need to make sure $\tilde{x}_t \in \X_{\epsilon_1}$ and $ \tilde{u}_t \in \U_{\epsilon_1}$. From \eqref{eq: tilde x u}, we equivalently need
\begin{equation}\label{eq: coupled safety constraint}
\begin{aligned}
    &\sum_{k=0}^{2H} \Psi_k^x (\M_{t-H:t-1}) w_{t-k} \in \X_{\epsilon_1},\\
    &\sum_{k=0}^{2H} \Psi_k^u (\M_{t-H:t}) w_{t-k} \in \U_{\epsilon_1}, \quad \forall w_{t-k} \in \W.
\end{aligned} 
\end{equation}

However, verifying the aforementioned safety constraints is challenging due to their coupling of all the policies \(\M_{t-H:t}\) from the past \(H\) steps. This issue can be addressed by approximating \(\Psi_k^x(\M_{t-H:t-1}), \Psi_k^u (\M_{t-H:t})\) with \(\mathring{\Psi}_k^x(\M_t), \mathring{\Psi}_k^u(\M_t)\) as defined in \eqref{psiring}. This approximation is reasonable if the policy \(\M_t\) changes slowly over time, which can be ensured by configuring our online algorithm to operate with a diminishing step size. Consequently, the error incurred by this alternative evaluation can be bounded by a small parameter $\epsilon_2$ provided that we select the step size to be sufficiently small (shown in Lemma \ref{lemma: epsilon_2 bound} in Appendix \ref{Appendix: thm1}). 
Therefore, we can ensure that the algorithm is safe by enforcing
$$
\sum_{k=0}^{2H} \mathring{\Psi}_k^x w_{t-k} \in \X_{\epsilon_1 + \epsilon_2}, \quad \sum_{k=0}^{2H} \mathring{\Psi}_k^u w_{t-k} \in \U_{\epsilon_1 + \epsilon_2},
$$
under any disturbance, which ensures that
$$
\begin{aligned}
    \tilde{x}_t &= \sum_{k=0}^{2H} \mathring{\Psi}_k^x w_{t-k} + \sum_{k=0}^{2H} (\Psi_k^x - \mathring{\Psi}_k^x) w_{t-k} \\
    &\in \X_{\epsilon_1 + \epsilon_2} + \ball^n(\epsilon_2) = (\X_{\epsilon_1 + \epsilon_2})_{-\epsilon_2} \subseteq \X_{\epsilon_1}.
\end{aligned}
$$
Similar analysis can also be applied to $\tilde{u}_t$, and we will have $\tilde{u}_t \in \U_{\epsilon_1}$ as expected. 

To help write the summations \eqref{eq: coupled safety constraint} in a compact form, we define the following matrices: 
\begin{equation} \label{equation: hx, hu}
\begin{aligned}
    h^x(\M_{t-H:t-1}) &= 
        \begin{bmatrix}
        \Psi_1^x  &\Psi_2^x  & ... &\Psi_{2H}^x 
        \end{bmatrix}, \\
    h^u(\M_{t-H:t}) &= 
    \begin{bmatrix}
    \Psi_1^u  &\Psi_2^u  & ... &\Psi_{2H}^u 
    \end{bmatrix}, \\
\end{aligned}
\end{equation}
and similarly write $\mathring{h}^x(\M), \mathring{h}^u(\M)$ as 
$$\mathring{h}^x(\M) =  h^x(\M,..., \M), \quad \mathring{h}^u(\M) =  h^u(\M,..., \M), $$
when a fixed weight matrix $\M$ is applied in the past $H$ step, where $h^x, \mathring{h}^x \in \R^{n \times 2Hn}$, and  $h^u, \mathring{h}^u \in \R^{m \times 2Hn}$.

Now, we are prepared to define our set of safe DAC policies. This safety set plays a critical role in our ability to select a DAC that consistently ensures the state and control input satisfy our safety concerns. 
The safety policy set $\Omega^{\epsilon}$ is defined as follows:
\begin{equation} \label{eq: Omega^epsilon}
\begin{aligned}
\Omega^{\epsilon} =  \big \{ \M \in \mathcal{M}: \quad &  \mathring{h}^x(\M) \cdot \ball (\overline{w}) \subseteq \X_{\epsilon} \\  
         &  \mathring{h}^u(\M) \cdot \ball (\overline{w}) \subseteq \U_{\epsilon}   \big \},
\end{aligned}
\end{equation}
where $\ball (\overline{w})$ represents the closed $L_{\infty}$ ball lying in $\R^{2Hn}$, corresponding to the size of $\mathring{h}^x$ and $\mathring{h}^u$, the set $\mathcal{M}$ is defined in \eqref{eq: M matrix}. The process of selecting an appropriate buffer size $\epsilon$ to ensure the algorithm's safety and low regret will be discussed in  Section \ref{section: thm}.

\textbf{Constructing the Approximate Cost Function:}

The method we use to construct an approximate cost function is based on the work in \cite{agarwal_online_2019}, where an online nonstochastic control problem is reformulated into an ``OCO with memory" problem, stated as follows: At each stage \(t\), the agent selects a policy \(\M_t \in \mathcal{M} \) and incurs a cost \(f_t(\M_{t-H:t})\). Despite the coupling of the current policy \(\M_t\) with historical policies \(\M_{t-H:t-1}\), $\mathcal{M}$ defined in \eqref{eq: M matrix} remains decoupled and depends solely on the current \(\M_t\).

To solve this ``OCO with memory" problem, the author of \cite{agarwal_online_2019} defines an approximate cost function $ \mathring f_t $ as:
    \begin{equation}
        f_t(\M_{t-H:t}) := c_t(\tilde{x}_t,\tilde{u}_t),\ 
        \mathring f_t(\M) := f_t(\M,..., \M). 
    \label{idealized cost}
    \end{equation}
    
Here, $f_t$ is called the idealized cost and is convex with respect to $\M_{t-H:t}$ since $\tilde{x}_t,\tilde{u}_t$ are affine functions of $\M_{t-H:t}$ and $c_t(\cdot,\cdot)$ is convex. For the same reason, the approximate cost function $\mathring f_t$ is also convex with respect to $\M$.
    
Consequently, the OCO with memory problem is reformulated as a classical OCO problem with stage cost \(\mathring{f}_t(\M_t)\), which can be tackled using classical OCO algorithms such as online gradient descent (OGD) \cite{zinkevich_online_2003}. The stepsize of OGD is chosen sufficiently small to ensure minimal variation between the current policy \(\M_t\) and the historical policies \(\M_{t-H}, \ldots, \M_t\), guaranteeing small approximation errors between \(\mathring{f}_t(\M_t)\) and \(f_t(\M_{t-H:t})\), and thus low regret.

With the groundwork laid in the preceding steps, we now proceed to formulate our online algorithm. 

\begin{algorithm2e}
	\caption{OGD-BZC}
	\label{alg:ogd}
	\SetAlgoNoLine
	\DontPrintSemicolon
	\LinesNumbered
	\KwIn{A $(\kappa,\gamma)$-strongly stable  matrix $K$, memory size $H>0$,  buffer size $\epsilon$,   stepsize $\eta$. }
	
	Initialize $\bm M_0 \in \Omega^\epsilon$ defined in \eqref{eq: Omega^epsilon}.\;

	\For{$t=0,  1,2,\ldots,T$}{
		Apply control input $	u_t=-K x_t \!+\! \sum_{i=1}^H \M_t^{[i]}w_{t-i}.
		$\;
		
		Observe the next state $x_{t+1}$ and calculate $w_t$.\;
		
		Run projected OGD 
		$$ \bm M_{t+1}=\Pi_{\Omega^\epsilon}\left[\bm M_t-\eta \nabla \mathring f_t(\bm M_t)\right],$$
		where $\mathring f_t(\bm M)$ is defined in \eqref{idealized cost}.}
	
\end{algorithm2e}

The parameters $\epsilon$, $H$, and $\eta$ must be carefully selected to guarantee small approximation errors in Steps 1-4 and achieve sub-linear regret. Conditions for suitable choices of parameters are discussed in  Section \ref{section: thm}.

\section{Theoretical Result}\label{section: thm}
In this section, we prove that OGD-BZC ensures safety requirements and achieve $\tilde{\mathcal{O}}(\sqrt{T})$ regret if the parameters are properly chosen. We will state our theoretical results in the following order. First, we show that OGD-BZC is safe in Theorem \ref{thm: safety}, and then provide a regret bound in Theorem \ref{thm: regret bound}.

\subsection{Safety Result}
In the following theorem, we prove that OGD-BZC is a safe algorithm with appropriate choice of parameters.

\begin{theorem} \label{thm: safety}
    Assume $\epsilon \ge 0$ and $H \ge \frac{\log(2\kappa^2)}{\log((1-\gamma)^{-1})}$. 
    Also, let
    $$
    \begin{aligned}
        \epsilon_1(H) & :=c_1 H (1-\gamma)^H,\\
        \epsilon_2(\eta,H) & := c_2 \sqrt{mn^3}\eta H^2,\\
        \epsilon_3(H) & := c_3 \sqrt{n} (1-\gamma)^H.
    \end{aligned}
    $$
    where $c_1 = \overline{w}\kappa^3(2\kappa^3 + 2a\kappa^3 \kappa_B + a)/\gamma$, $c_2 = 4 G \overline{w}^3 (\kappa^3 + 2a\kappa^3 \kappa_B)(1+\kappa)\kappa^5\kappa_B^2/\gamma^4$ and $c_3 = 2\overline{w}\kappa^5$.
    Then, if the buffer size $\epsilon$ and $H$ meet the condition
        $$ 
        \begin{aligned}
            \epsilon_1(H) + \epsilon_2(\eta,H) \le \epsilon \le \epsilon_* - \epsilon_1(H) - \epsilon_3(H),
        \end{aligned}  
        $$
    then the safety policy set $\Omega^{\epsilon}$ is non-empty, and OGD-BZC algorithm is safe under any sequence of disturbances $\left \{ w_k \in \W \right \}_{k=0}^T$.
\end{theorem}
The proof of Theorem \ref{thm: safety} differs from \cite{li_online_2021} primarily in how we define the notion of ``buffer zone" based on set shrinkage and expansion to handle convex sets. We will see that this distinct approach introduces additional challenges in the regret analysis presented in Section \ref{subsec: regret}.

\subsection{Regret Analysis}\label{subsec: regret}
In this subsection, we analyze the regret of the OGD-BZC algorithm and detail the specific parameters needed for the algorithm to reach a regret of $\tilde{\mathcal{O}}(\sqrt{T})$. 

\begin{theorem}[Regret Bound] \label{thm: regret bound}
    Under the requirement of $H$ and buffer size $\epsilon$ in Theorem \ref{thm: safety}, if the buffer size additionally satisfies 
    $$ \epsilon < \frac{\epsilon_*}{2} - \epsilon_1(H) - \epsilon_3(H),$$
    OGD-BZC will have the following regret bound:

    $$
    \begin{aligned}
        \text{Reg}_T (\text{OGD-BZC}) \le \mathcal{O} \bigg ( & T \sqrt{m^2 n^3} H^2 (1-\gamma)^H + Tmn^2H^3\eta  \\
        & + \frac{m}{\eta} + T(1-\gamma)^H \sqrt{m^2n^3H^5}/\epsilon_* \\
        & + \epsilon T\sqrt{m^2n^2H^3}/\epsilon_* \bigg ),
    \end{aligned}
    $$
where the hidden constant coefficients is the polynomial of $\overline{w}, G, \kappa, \kappa_B, a, \gamma^{-1}$.
\end{theorem}

\begin{corollary}[Sublinear Regret]
For sufficiently large $T$, if we set
$$H = \frac{\log ((8c_1T+ 4c_3\sqrt{n})/\epsilon_*)}{\log ((1-\gamma)^{-1})} ,$$ 
and $\eta = \frac{1}{n\sqrt{T H^3}}, 
\epsilon = \epsilon_1(H) + \epsilon_2 (\eta, H)$, then OGD-BZC is safe with
$$ \text{Reg}_T(\text{OGD-BZC}) \le \tilde{\mathcal{O}} (m^{1.5} n^{1.5} \sqrt{T}). $$
\end{corollary}

The regret result derived here has the same order as the result in the unconstrained online nonstochastic control literature for convex cost functions \cite{agarwal_online_2019}. 

\subsection{Proof of Theorem \ref{thm: regret bound}} 
Let $J_T(\A)$ be the total loss received by our algorithm, i.e.  $J_T(\A) = \sum_{t=0}^T c_t\left(x_t^{\A}, u_t^{\A}\right)$, then the proof relies on dividing the regret into two parts and give bound guarantee for each part under any sequence of disturbances. Specifically,
$$
\begin{aligned} 
& J_T(\A)-\min_{K \in \K} J_T(K) =  \\
& \underbrace{J_T(\A)- \min _{\M \in \Omega^{\epsilon}} \sum_{t=0}^T \mathring{f}_t(\M) }_{\text{Performance Gap}} + \underbrace{\min _{\M \in \Omega^{\epsilon}} \sum_{t=0}^T \mathring{f}_t(\M)- \min_{K \in \K} J_T(K)}_{\text{Policy Set Gap}}.
\end{aligned}
$$

We will first provide a bound for the ``Policy Set Gap" part. In order to do so, we define the DAC approximation for a given linear policy $K' \in \K$ as
$$
\M^{[i]}(K'):= (K -K')(A - BK')^{i-1}.
$$

Also, we use the notation $\M_{\text{ap}} = \M(K^*)$ where $K^* := \text{argmin}_{K \in \K} J_T(K)$.
Then the ``Policy Set Gap'' can be separated as
$$
\begin{aligned}
    &\min _{\M \in \Omega^{\epsilon}} \sum_{t=0}^T \mathring{f}_t(\M)- \min_{K \in \K} J_T(K) =    \\ 
    & \underbrace{\min _{\M \in \Omega^{\epsilon}} \sum_{t=0}^T \mathring{f}_t(\M) -  \mathring{f}_t(\M_{\text{ap}})}_{(\bigstar)}  + \underbrace{\sum_{t=0}^T \mathring{f}_t(\M_{\text{ap}})  - \min_{K \in \K} J_T(K)}_{\text{Lemma \ref{lemma: M_ap}}}.
\end{aligned}
$$

An error bound for the second term is given in Lemma \ref{lemma: M_ap}, which will be discussed later in this section. Now it remains to bound the first term, denoted as $(\bigstar)$. 

We highlight that bounding $(\bigstar)$ is one of the main contributions of this paper, as the method to derive the bound here is distinct from existing works and is highly nontrivial. Specifically, in \cite{li_online_2021}, the authors derive a bound in the affine context by leveraging the specific structure of affine constraints. This approach results in a bound that scales with the number of affine constraints, which cannot be directly applied in the convex setting.

Here, we propose Lemma \ref{lemma: star} that (i) provides a regret bound that applies to general convex constraints and (ii) removes the dependence of the regret bound on the number of constraints.
\begin{lemma}\label{lemma: star}
Under the condition of Theorem \ref{thm: regret bound}, we have
    $$ (\bigstar) \le \mathcal{O} \left ( T \sqrt{m^2 n^2 H^3} \frac{(\epsilon_1 + \epsilon_3 + \epsilon)}{\epsilon_*} \right),$$
under any sequence of disturbances.
\end{lemma}

To prove Lemma \ref{lemma: star}, we depend on an auxiliary lemma that estimates the size difference between $\Omega^{\epsilon}$ and $\Omega^{-\epsilon_1 - \epsilon_3}$, denoted as Lemma \ref{lemma: Omega set bound}. This estimation is achieved by utilizing tools from classical convex analysis and the proof can be found in the Appendix \ref{Appendix: thm2}.
\begin{lemma} \label{lemma: Omega set bound}
    If the buffer size $\epsilon$ and $H$ satisfies
    $$ \epsilon < \epsilon_* - \epsilon_1(H) - \epsilon_3(H),$$
    there exists $\alpha_1 = 1 - \frac{\epsilon}{\epsilon+r}, \quad \alpha_2 = 1 + \frac{\epsilon_1 + \epsilon_3}{\epsilon+r}$ such that 
        $$ \mathring \Omega(\alpha_1) \subseteq \Omega^{\epsilon}, \quad  \mathring \Omega (\alpha_2)\supseteq \Omega^{-\epsilon_1 -\epsilon_3},$$
    where $r :=  \epsilon_* - (\epsilon + \epsilon_1 + \epsilon_3)$, and $ \mathring \Omega(\alpha) := \left \{ \alpha \M + (1-\alpha)\M(K_{ss}): \M \in \Omega   \right\}$.

\end{lemma}

\begin{proof}[Proof of Lemma \ref{lemma: star}]


To bound $(\bigstar)$, we first need to utilize the fact that $\mathring{f}_t$ is convex and possesses a bounded gradient. The property of $\mathring{f}_t$ having a bounded gradient originates from \cite{agarwal_online_2019}. In the context of this paper, we demonstrate that the bound for $||\nabla \mathring{f}_t||_F$ is $G_f = \mathcal{O}(\sqrt{mn^2H^3})$, where $\norm{\cdot}_F$ denotes the Frobenius norm (see Lemma \ref{lemma: bound on f_t} in Appendix \ref{Appendix: helping}) . Then by the bounded gradient and the convexity of $\mathring{f}_t$, we have $(\bigstar) \le \min _{\M \in \Omega^{\epsilon}} T G_f \norm{\M - \M_{\text{ap}}}_F$.

Next, we will use Lemma \ref{lemma: Omega set bound} to bound the distance between $\M_{\text{ap}}$ and $\Omega^\epsilon$.
In particular, we have that $\M_{\text{ap}} \in \Omega^{-\epsilon_1 - \epsilon_3} \subseteq \mathring \Omega (\alpha_2)$, and therefore, there exists $\Tilde{\M} \in \Omega$ such that 
\begin{equation}
\begin{split}
\label{eqn:map}
    & \M_{\text{ap}} = \alpha_2 \Tilde{\M} + (1-\alpha_2)\M(K_{ss}).
\end{split}
\end{equation}

At the same time, $\mathring \Omega(\alpha_1) \subseteq \Omega^{\epsilon}$, and therefore 
\begin{equation}
\label{eqn:mtil}
    \alpha_1 \Tilde{\M} + (1-\alpha_1)\M(K_{ss}) \in \Omega^{\epsilon}.
\end{equation}

Therefore, combining \eqref{eqn:map} and \eqref{eqn:mtil}, we have that
\begin{equation*}
    \frac{\alpha_1}{\alpha_2} \M_{\text{ap}} + \left(1-\frac{\alpha_1}{\alpha_2} \right)\M(K_{ss})  \in \Omega^{\epsilon}.
\end{equation*}

It follows that
$$\begin{aligned}
    &\min _{\M \in \Omega^{\epsilon}} T G_f \norm{\M - \M_{\text{ap}}}_F \\
    & \le T G_f \norm{\frac{\alpha_1}{\alpha_2}\M_{\text{ap}} + \left(1 - \frac{\alpha_1}{\alpha_2} \right) \M(K_{ss}) - \M_{\text{ap}}}_F \\
    &= T G_f \left(1 - \frac{\alpha_1}{\alpha_2} \right)\norm{\M_{\text{ap}} - \M(K_{ss})}_F \\
    & \le T G_f \left(\frac{\epsilon + \epsilon_1 + \epsilon_3}{r}\right) \norm{\M_{\text{ap}} - \M(K_{ss})}_F .
\end{aligned}
$$

To bound $\norm{\M_{\text{ap}} - \M(K_{ss})}_F$, we use the fact that $\M_{\text{ap}}$, $\M(K_{ss}) \in \mathcal{M}$. As a result,
$$
\begin{aligned}
&\norm{\M_{\text{ap}} - \M(K_{ss})}^2_F = \sum_{i=1}^H \norm{\M_{\text{ap}}^{[i]} - \M^{[i]}(K_{ss})}^2_F \\
&\le \sum_{i=1}^H m \norm{\M_{\text{ap}}^{[i]} - \M^{[i]}(K_{ss})}^2_2 \\
&\le \sum_{i=1}^H m \left (\norm{\M_{\text{ap}}^{[i]}}_2 + \norm{\M^{[i]}(K_{ss})}_2 \right )^2 \\
&\le \sum_{i=1}^H m \cdot (2a(1-\gamma)^{i-1})^2 \le 4ma^2/\gamma.
\end{aligned}
$$

Putting everything together, and using the fact that $r = \epsilon_* - (\epsilon + \epsilon_1 + \epsilon_3) > \epsilon_*/2 $, it follows that
$$
\begin{aligned}
    (\bigstar) & \leq T G_f \left(\frac{\epsilon + \epsilon_1 + \epsilon_3}{r}\right) \sqrt{4ma^2/\gamma}\\
    & < 2 T G_f \left(\frac{\epsilon + \epsilon_1 + \epsilon_3}{\epsilon_*}\right) \sqrt{4ma^2/\gamma}.
\end{aligned}
$$

Then, since $G_f = \mathcal{O}(\sqrt{mn^2H^3})$, we conclude that $(\bigstar)$ is bounded by $ \mathcal{O}\left ( T \sqrt{m^2 n^2 H^3} (\epsilon_1 + \epsilon_3 + \epsilon)/\epsilon_* \right  )$.
\end{proof}

\begin{lemma} \label{lemma: M_ap}
Given $K^* \in \K$ and $\M_{\text{ap}} = \M (K^*)$, we will have 
$\M_{\text{ap}} \in \Omega^{-\epsilon_1 - \epsilon_3}$ and 
$$
\sum_{t=0}^T \mathring{f}_t(\M_{\text{ap}})  - \min_{K \in \K} J_T(K) \le \mathcal{O} \left( TmnH (\epsilon_1 + \epsilon_3) \right ),
$$
under any sequence of disturbances.
\end{lemma}

 The proof of Lemma \ref{lemma: M_ap} is based on \cite{li_online_2021} with minor differences in coefficient and notation from its original proof. The comprehensive proof of Lemma \ref{lemma: M_ap} is provided in the Appendix~\ref{Appendix: thm2}.

To finish the proof of Theorem \ref{thm: regret bound}, we now derive a bound the ``Performance Gap'', as provided in the following lemma.

\begin{lemma}[Performance Gap]\label{lemma: Performance gap}
   Under the condition of Theorem \ref{thm: regret bound}, we have
    $$
    \begin{aligned}
    J_T(\A)- \min _{\M \in \Omega^{\epsilon}} \sum_{t=0}^T \mathring{f}_t(\M)
    \le \mathcal{O} \Big (  & T mn H^2 (1-\gamma)^H \\
    & + Tmn^2H^3\eta + \frac{m}{\eta} \Big),
    \end{aligned}
    $$
    under any sequence of disturbances.
\end{lemma}
The proof of Lemma \ref{lemma: Performance gap} relies on bounding two separated parts. The first one is to bound the gap between the true cost $c_t\left(x_t^{\A}, u_t^{\A}\right)$ and the approximated cost $\mathring f_t (\M_t)$. This bound can be found in  \cite{agarwal_online_2019} as direct result from the unconstrained online nonstochastic control setting. The second part involves estimating the remaining part by performing classical OGD analysis in the literature \cite{hazan_introduction_2023}, with $\mathring f_t$ serving as the objective function. The detailed proof of Lemma \ref{lemma: Performance gap} can be found in the Appendix \ref{Appendix: thm2}, with subtle notation difference from their original proof in \cite{agarwal_online_2019, hazan_introduction_2023, li_online_2021}.

Finally, by combining the results from Lemma \ref{lemma: star}, Lemma \ref{lemma: M_ap}, and Lemma \ref{lemma: Performance gap}, we establish that the regret of OGD-BZC is of the form as stated in Theorem \ref{thm: regret bound}.

\section{Numerical Experiments} \label{section: numerical}

In this section, we present the numerical experiments conducted to evaluate the performance of our proposed algorithm. The experiments aim to demonstrate the efficacy and robustness of the algorithm in addressing constrained online non-stochastic control problems.

We will apply the OGD-BZC algorithm in a toy example with 2-D dynamics described as follows:
$$
    x_{t+1} = \begin{bmatrix}
1 & 1 \\
0 & 1/2 \\
\end{bmatrix}x_t 
+ \begin{bmatrix}
1  \\
1  \\
\end{bmatrix} u_t + w_t, \quad \forall t \ge 0,
$$
where \(x_t \in \mathbb{R}^2\) and \(u_t \in \mathbb{R}\). 
For our initial conditions, we set \(x_0 = u_0 = 0\). The safety constraints is set to be \(\|x_t\|_2 \leq 1\) and \(\|u_t\|_2 \leq 1\). We consider adversarial disturbances \(w_t \in \mathcal{W} := \{w: \|w\|_{\infty} \leq \overline{w}\}\) with $\overline{w} = 0.3$.  
For the algorithm's parameters, the matrix \(K = \begin{bmatrix}
1/2 & 0
\end{bmatrix}\) is chosen to be strongly stable matrix. Then, we set the memory size \(H = \lfloor \log(T) \rfloor \) and the buffer size \(\epsilon = \log(T)/\sqrt{T}\). The stepsize is defined as \(\eta = 1/(\sqrt{T}\log(T))\). We also consider the cost function $c_t(x,u) = \norm{x}_2^2 + \norm{u}_2^2$ for all $0 \le t \le T$.

\begin{figure}
\centering
    \subfigure[Random disturbances]{\includegraphics[width=0.45\linewidth]{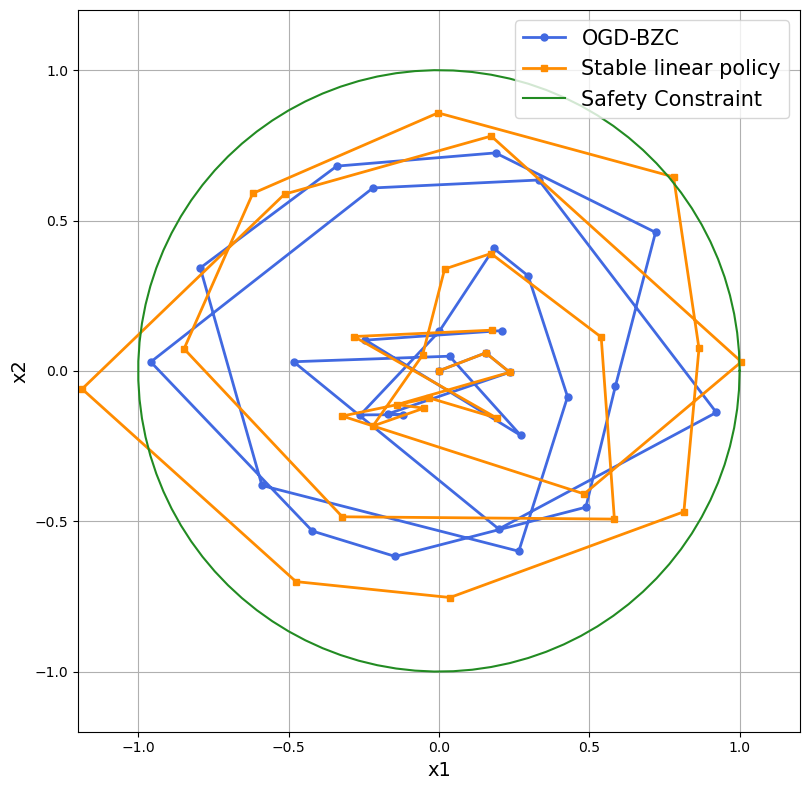}}\quad
    \subfigure[Assigned disturbances]{\includegraphics[width=0.45\linewidth]{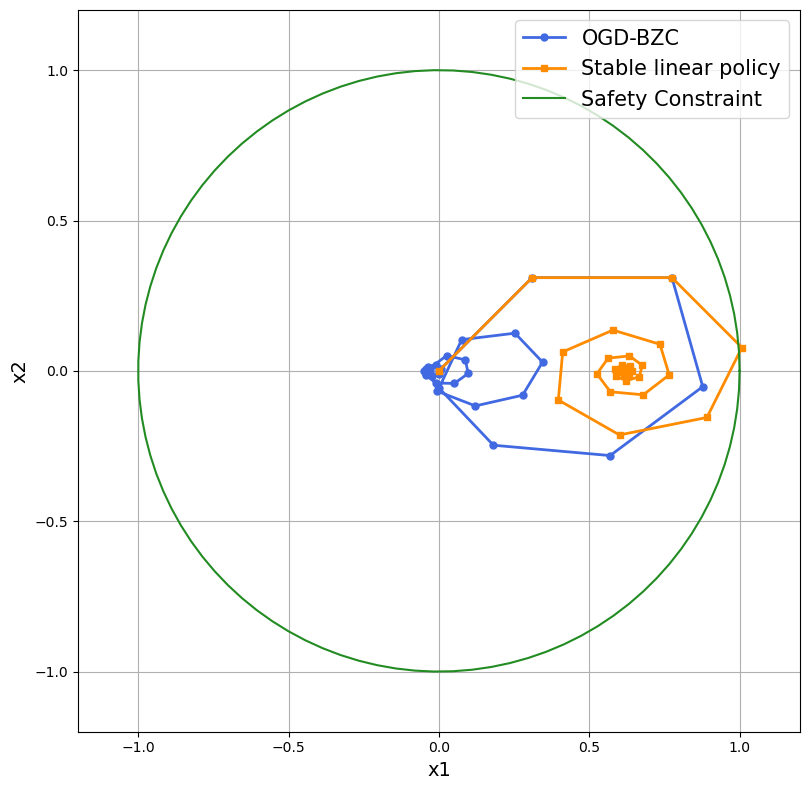}}
\caption{State Evolution of OGD-BZC and Stable linear controller with horizon $T = 30$.}
    \vspace{-0.4cm}
\label{fig:evolution}
\end{figure}

\begin{figure}
    \centering
    \includegraphics[width=0.75\linewidth]{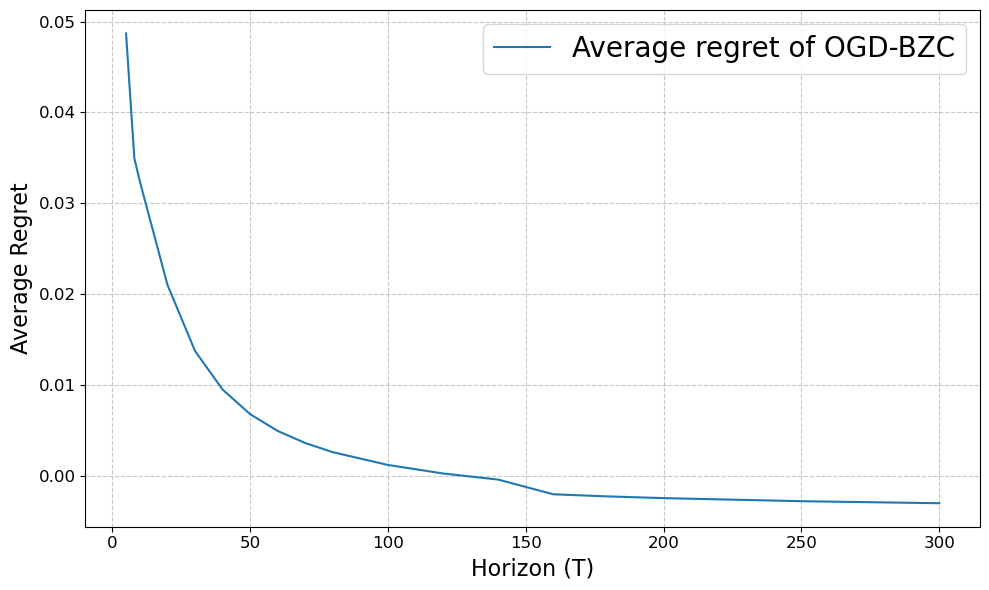}
    \caption{Average regret of OGD-BZC.}
    \label{fig:regret}
    \vspace{-0.4cm}
\end{figure}

The safety results are shown in Figure \ref{fig:evolution}, while the regret performance is illustrated in Figure \ref{fig:regret}.
Figure \ref{fig:evolution} depicts the evolution of the state under OGD-BZC and under a strongly stable linear controller (which is simply chosen to be $K$), respectively. Two scenarios are presented: one with i.i.d disturbances generated uniformly from $[-\overline{w},\overline{w}]$ and the other with assigned disturbances $w_t = \begin{bmatrix}
\overline{w} & \overline{w}
\end{bmatrix}^T, \forall 0 \le t \le T$.
In Figure \ref{fig:evolution}(a), it is observed that the trajectory of states generated by OGD-BZC remains within the prescribed safety set throughout the time horizon. 
In Figure \ref{fig:evolution}(b), in addition to the safety constraints being satisfied for all $0 \le t \le T$, it is also observed that OGD-BZC ``learns" from the assigned disturbance sequence and converges to points closer to the origin compared to the strongly stable linear policy.
Figure \ref{fig:regret} displays the average regret $\text{Reg}_T/T$ of OGD-BZC compared with the best safe linear policy in hindsight. The disturbance for calculating regret is set as $w_t = \begin{bmatrix} \overline{w} & \overline{w} \end{bmatrix}^T, \forall 0 \le t \le T$, to simplify the comparison. Notably, in our example, the regret indeed converges to a value smaller than $0$. This occurrence is likely not due to numerical error but because DAC encompasses a larger policy set than the linear policy set to which it is compared. Thus, it's feasible that a DAC outperforms the best safe linear policy in hindsight.

\section{Conclusion} \label{section: conclusion}
This paper studies the problem of online nonstochastic control for LTI systems under general convex state and control input constraints. Through the development and analysis of the OGD-BZC algorithm, we have generalized the work established by \cite{li_online_2021} to handle broader safety constraints. Our work demonstrates that, with carefully chosen parameters, OGD-BZC can ensure convex safety constraints while achieving the $\tilde{\mathcal{O}}(\sqrt{T})$ regret. Our numerical experiments reinforce the theoretical findings, showcasing OGD-BZC's robustness and effectiveness in practical scenarios.

Looking ahead, there are several directions to be explored or improved in the future. For instance, we can conduct further research on how to decrease the computational burden in the projection step of our algorithm. We can also consider using other policy sets instead of DAC. As an exploration direction, we can consider safe online nonstochastic control with static but unknown constraints and limited feedback at each step on the constraints.

\section{Acknowledgements}
This work was supported by NSF grant \#2330154.

\bibliographystyle{IEEEtran}
\bibliography{IEEEabrv,references}

\newpage
\onecolumn

\appendix

\subsection{Supporting Lemmas} \label{Appendix: support}
    



\begin{proposition} \label{proposition: Minkwski subtraction}
    Let $\D$ be a closed convex set in $\R^d$ and $\norm{\cdot}$ be any vector norm. Then, for any $\Delta_1, \Delta_2 > 0$, the sets $\D_{\Delta_1}$ and $\D_{-\Delta_1}$ are closed convex, and the following relations hold:
    \begin{equation}
    \begin{aligned}
    (\D_{\Delta_1})_{\Delta_2} &= \D_{\Delta_1 + \Delta_2}, \\
    (\D_{-\Delta_1})_{-\Delta_2} &= \D_{-(\Delta_1 + \Delta_2)}, \\
    (\D_{-\Delta_1})_{\Delta_2} &= \D_{\Delta_2 - \Delta_1}, \\
    (\D_{\Delta_2})_{-\Delta_1} &\subseteq \D_{\Delta_2 - \Delta_1}. \\
    \end{aligned}
    \end{equation}
\end{proposition}
    
 As a simple corollary, if we set $\Delta_1 = \Delta_2 = \Delta $ , the following equality holds for convex sets: $(\D_{\Delta})_{-\Delta} \subseteq \D = (\D_{-\Delta})_{\Delta} $. It implies that if when we first shrink some safety constraint $\D$ to $\D_{\Delta}$, and then allow an expansion with radius less or equal to $\Delta$, eventually we will still stay in the safety constraint. 
\begin{proof}
Firstly, if $\D$ is closed and convex, then $\D_{\Delta_1}$ is an intersection of closed and convex sets and hence is also closed and convex. It is also direct to see that $(\D_{\Delta_1})_{\Delta_2} = \D_{\Delta_1 + \Delta_2}$, and $(\D_{-\Delta_1})_{-\Delta_2} = \D_{-(\Delta_1 + \Delta_2)}$ by simply following the definition.

Next, we prove that $(\D_{\Delta_2})_{-\Delta_1} \subseteq \D_{\Delta_2 - \Delta_1}$.
For any $x\in \D_{\Delta_2}$, by definition we have $x + \ball(\Delta_2) \subseteq \D $ and 
$ x + \ball (\Delta_1) \subseteq (\D_{\Delta_2})_{-\Delta_1}$. If $\Delta_1 \ge \Delta_2$, then 
$$ \begin{aligned}
x + \ball(\Delta_1) &= x + \ball(\Delta_2) + \ball(\Delta_1 - \Delta_2) \\
&\subseteq D + \ball(\Delta_1 - \Delta_2) \\
&= \D_{\Delta_2 - \Delta_1}. 
\end{aligned}
$$

This proves $(\D_{\Delta_2})_{-\Delta_1} \subseteq \D_{\Delta_2 - \Delta_1}$. If $\Delta_1 < \Delta_2$, then
$$ \begin{aligned}
& x + \ball(\Delta_2) \subseteq \D \\
\Leftrightarrow \quad & x + \ball(\Delta_1) + \ball(\Delta_2 - \Delta_1) \subseteq \D \\
\Leftrightarrow \quad & x + \ball(\Delta_1) \subseteq \D_{\Delta_2-\Delta_1}.
\end{aligned}
$$

In conclusion, we have $(\D_{\Delta_2})_{-\Delta_1} \subseteq \D_{\Delta_2 - \Delta_1}$. 

Lastly, we prove that $(\D_{-\Delta_1})_{\Delta_2} = \D_{\Delta_2 - \Delta_1}$. 
If $\Delta_1 \ge \Delta_2$, we have 
$$ \begin{aligned}
&\D + \ball(\Delta_1) = \D + \ball(\Delta_1 - \Delta_2) + \ball(\Delta_2) = \D_{-\Delta_1} \\
\Rightarrow \quad & \D + \ball(\Delta_1 - \Delta_2) \subseteq (\D_{-\Delta_1})_{\Delta_2}.
\end{aligned}
$$

This implies $(\D_{-\Delta_1})_{\Delta_2} \supseteq \D_{\Delta_2 - \Delta_1}$. If $\Delta_1 < \Delta_2$, for any $x \in \D_{\Delta_2 - \Delta_1} $, we have 
$$ \begin{aligned}
& x + \ball(\Delta_2 - \Delta_1) \subseteq \D \\
\Rightarrow \quad & x + \ball(\Delta_2) = x + \ball(\Delta_2 - \Delta_1) + \ball(\Delta_1) \subseteq \D_{-\Delta_1} \\
\Rightarrow \quad & x \in (\D_{-\Delta_1})_{\Delta_2}. \\
\end{aligned}
$$

So far we already proved $(\D_{-\Delta_1})_{\Delta_2} \supseteq \D_{\Delta_2 - \Delta_1}$. For the reverse direction, we consider the \textit{support function} of closed and convex set:
$$
h_{\D}(y) := \sup \{ \langle y,u \rangle: \forall u \in \D \} \quad \forall y \in \R^d.
$$

For $x\in (\D_{-\Delta_1})_{\Delta_2}$, we have $x + \ball(\Delta_2) \subseteq D + \ball(\Delta_1) $. By the property supporting function (one can see \cite{schneider_convex_2014} for more details), we have $h_{\{ x\}} + h_{\ball(\Delta_2)} \le h_{\D} + h_{\ball(\Delta_1)}$. Also,it is easy to verify that the support function of a ball is a scaled norm function. (i.e. $h_{\ball(\Delta)}(y) = \Delta \norm{y}$). Thus we have
$$
    h_{\{ x\}}(y) + \Delta_2 \norm{y} \le h_{\D} + \Delta_1 \norm{y}.
$$

If $\Delta_1 \ge \Delta_2$, then $ h_{\{ x\}}(y) \le h_{\D} + (\Delta_1-\Delta_2) \norm{y}$ implies $x \in \D + \ball(\Delta_1 - \Delta_2)$, which is $x \in \D_{\Delta_2 - \Delta_1}$. 

If $\Delta_1 < \Delta_2$, then $ h_{\{ x\}}(y) + (\Delta_2-\Delta_1) \norm{y}\le h_{\D} $ implies implies $x + \ball(\Delta_2 - \Delta_1) \in \D $, which is also $x \in \D_{\Delta_2 - \Delta_1}$.
\end{proof}

\newpage

\subsection{Proofs related to Theorem \ref{thm: safety}}
\label{Appendix: thm1}

Before giving Theorem \ref{thm: safety}, we give three lemmas that will be used by Theorem \ref{thm: safety}.

\begin{lemma}\label{lemma: epsilon_1 bound}
    Let $x_t, u_t$ represents the state and control input generated by OGD-BZC, while $\tilde{x}_t, \tilde{u}_t$ denote the surrogate state and input defined in \eqref{prop: DAC}. 
    Then,
    under the condition of Theorem \ref{thm: safety}, it holds for any disturbance sequence $\{w_k \in \mathcal{W}\}$ that 
    $$
    \begin{aligned}
    & x_t - \tilde{x}_t \in \ball (\epsilon_1(H)), \\
    & u_t - \tilde{u}_t \in \ball (\epsilon_1(H)),
    \end{aligned}
    $$
\end{lemma}
\begin{proof}
By Lemma \ref{lemma: x_t} we obtain that $\max \{\norm{x_t}_2, \norm{u_t}_2 \} \le b$. Thus
$$
\begin{aligned}
    \norm{x_t - \tilde{x}_t}_{\infty} = \norm{A_K^H x_{t-H}}_{\infty} &\le \norm{A_K^H x_{t-H}}_2 \le \norm{A_K^H}_2 \norm{x_{t-H}}_2 \\
    & \le \kappa^2 (1-\gamma)^H b. \\
    \norm{u_t - \tilde{u}_t}_{\infty} = \norm{KA_K^H x_{t-H}}_{\infty} &\le \norm{KA_K^H x_{t-H}}_2 \le \norm{K}_2 \norm{A_K^H}_2 \norm{x_{t-H}}_2 \\
    & \le \kappa^3 (1-\gamma)^H b.
\end{aligned}
$$

Since when $H \ge \frac{\log (2\kappa^2)}{\log((1-\gamma)^{-1})}$, by Lemma \ref{lemma: x_t} we have
$$b \le 2 \overline{w}\sqrt{n} (\kappa^3 + a\kappa^3 \kappa_B H)/\gamma + \overline{w}\sqrt{m}  \cdot a /\gamma \le \mathcal{O}(\sqrt{mn}H). $$

We can set $\epsilon_1(H) = c_1 \sqrt{mn}(1-\gamma)^H H$, where $c_1 = 2 \overline{w}\kappa^3 (\kappa^3 + a\kappa^3 \kappa_B + a/2 )/\gamma$. Then $\norm{x_t - \tilde{x}_t}_{\infty} \le \epsilon_1(H)$ and $\norm{u_t - \tilde{u}_t}_{\infty} \le \epsilon_1(H)$.
It trivially follows that $x_t - \tilde{x}_t \in \ball (\epsilon_1(H))$ and $u_t - \tilde{u}_t \in \ball (\epsilon_1(H))$, completing the proof.
\end{proof}

\begin{lemma}\label{lemma: epsilon_2 bound}
Let $\M_t$ represents the weight matrix generated by OGD-BZC. 
Then, under the condition of Theorem \ref{thm: safety}, we have
    $$\begin{aligned}
        (h^x(\M_{t-H+1:t}) - \mathring{h}^x(\M_t)) \cdot \ball^{2Hn}(\overline{w}) & \subseteq \ball^{n} (\epsilon_2 (\eta, H)),\\
        (h^u(\M_{t-H:t}) - \mathring{h}^u(\M_t)) \cdot \ball^{2Hn}(\overline{w}) & \subseteq \ball^{n} (\epsilon_2 (\eta, H)).\\
    \end{aligned}
    $$
\end{lemma}
\begin{proof}
By the definition of $h^x$ and $\mathring{h}^x$ in \eqref{equation: hx, hu} and Lemma \ref{lemma: Psi difference}, we have 
    \begin{equation}\label{eqn:hs3}\begin{aligned}
        \norm{h^x(\M_{t-H+1:t}) - \mathring{h}^x(\M_t)}_{\infty} &\le \sum_{k=1}^{2H} 
        \norm{\Psi_k^x(\M_{t-H+1:t}) - \mathring{\Psi}_k^x(\M_t)}_{\infty}\\
        & \le \sum_{k=1}^{2H} \sqrt{n} \norm{\Psi_k^x(\M_{t-H+1:t}) - \mathring{\Psi}_k^x(\M_t)}_2\\
        & \le \kappa^2 \kappa_B \sqrt{nH} \sum_{i=1}^{H}  (1 - \gamma)^{i-1} \norm{ \M_t - \M_{t-i} }_F
    \end{aligned}
    \end{equation}
Then, we bound the last term as
\begin{equation}\label{eqn:Ms}
\begin{aligned}
\norm{\M_t- \M_{t-i}}_F & \le \sum_{s=1}^{i} \norm{\M_{t-s+1}- \M_{t-s}}_F\\
& = \sum_{s=1}^{i} \norm{\Pi_{\Omega^\epsilon}\left[\bm M_{t-s}-\eta_{t-s} \nabla \mathring f_t(\bm M_{t-s})\right]- \M_{t-s}}_F\\
& \leq \sum_{s=1}^{i} \norm{\bm M_{t-s}-\eta_{t-s} \nabla \mathring f_t(\bm M_{t-s})- \M_{t-s}}_F\\
& \le 
\sum_{s=1}^{i} \eta \norm{\nabla \mathring{f}_{t-s}(\M_{t-s})}_F\\
& \le i\eta G_f,
\end{aligned}
\end{equation}
where the equality uses the update from the algorithm, the second inequality follows from the fact that $\Omega^\epsilon$ is convex,\footnote{To see that $\Omega^\epsilon$ is convex, first note that the set $\{ h\in \mathbb{R}^{n \times 2 H n} :  h \cdot \ball (\overline{w}) \subseteq \X_{\epsilon} \}$ is convex because it is the intersection of convex sets (given that $\X_{\epsilon}$ is convex). Therefore, by an affine mapping, the set $\{ \M \in \mathcal{M}: \mathring{h}^x(\M) \cdot \ball (\overline{w}) \subseteq \X_{\epsilon} \}$ is also convex. Since this also holds for the constraint on $\mathring{h}^u$, the set $\Omega^\epsilon$ is convex.} the third inequality holds because the step size is fixed (i.e. $\eta_{t-s} = \eta$) and the fourth inequality is from Lemma \ref{lemma: bound on f_t}.
Combining \eqref{eqn:hs3} and \eqref{eqn:Ms}, we have
    $$\begin{aligned}
        \norm{h^x(\M_{t-H+1:t}) - \mathring{h}^x(\M_t)}_{\infty} & \le \kappa^2 \kappa_B \sqrt{nH} \sum_{i=1}^{H}  (1 - \gamma)^{i-1} i\eta G_f \\
        & \le \kappa^2 \kappa_B \sqrt{nH} \eta G_f \frac{1}{\gamma^2}\\
        & \le \mathcal{O} (\sqrt{nH} \eta \cdot \sqrt{mn^2H^3}) \\
        & = \mathcal{O} (\sqrt{mn^3}H^2 \eta),
    \end{aligned}
    $$
    where we use the fact that $Gf \le \mathcal{O}(\sqrt{mn^2H^3})$.

Similarily, we can bound the distance between $h^u$ and $\mathring{h}^u$, as
    $$\begin{aligned}
        \norm{h^u(\M_{t-H:t}) - \mathring{h}^u(\M_t)}_{\infty} & \le \sum_{k=1}^{2H} 
        \norm{\Psi_k^u(\M_{t-H:t}) - \mathring{\Psi}_k^u(\M_t)}_{\infty}\\
        & \le \sum_{k=1}^{2H} \sqrt{m} \norm{\Psi_k^u(\M_{t-H:t}) - \mathring{\Psi}_k^u(\M_t)}_2\\
        & \le \kappa^3 \kappa_B \sqrt{mH} \sum_{i=1}^H (1-\gamma)^{i-1} \norm{ \M_t - \M_{t-i} }_F \\
        & \le \kappa^3 \kappa_B \sqrt{mH} \sum_{i=1}^{H} (1 - \gamma)^{i-1} i\eta G_f \\
        & \le \kappa^3 \kappa_B \sqrt{mH} \eta G_f \frac{1}{\gamma^2}\\
        & \le \mathcal{O} (\sqrt{mH} \eta \cdot \sqrt{mn^2H^3}) \\
        & = \mathcal{O} (mnH^2 \eta).
    \end{aligned}
    $$

Since $\epsilon_2(\eta, H) := \overline{w} \kappa^3 \kappa_B \sqrt{mnH} \eta G_f \frac{1}{\gamma^2} = c_2 \sqrt{m^2n^3}H^2 \eta$, it follows that 
$$
\norm{h^u(\M_{t-H:t}) - \mathring{h}^u(\M_t)}_{\infty} \overline{w} \leq \epsilon_2(\eta, H) \quad \text{and} \quad \norm{h^x(\M_{t-H+1:t}) - \mathring{h}^x(\M_t)}_{\infty} \overline{w} \leq \epsilon_2(\eta, H).
$$
Then, by the definition of the infinity matrix-norm, it follows that
$$
(h^x(\M_{t-H+1:t}) - \mathring{h}^x(\M_t)) \cdot \ball^{2Hn}(\overline{w}) \subseteq \ball^{n} (\epsilon_2 (\eta, H)), \quad \text{and} \quad
(h^u(\M_{t-H:t}) - \mathring{h}^u(\M_t)) \cdot \ball^{2Hn}(\overline{w}) \subseteq \ball^{n} (\epsilon_2 (\eta, H)),
$$
completing the proof.

\end{proof}

\begin{lemma}\label{lemma: epsilon_3 bound}
Given some $K' \in \K$, let
$$
\M^{[i]}(K'):= (K -K')(A - BK')^{i-1},
$$
and $x_t^{\M(K')},u_t^{\M(K')}$ be the states and inputs generated with disturbance-action policy $\M(K')$.
Then, for any $K' \in \K$, it holds that
    $$
    \begin{aligned}
    &x_t^{K'} - x_t^{\M(K')} \in \ball (\epsilon_3(H)), \\
    &u_t^{K'} - u_t^{\M(K')} \in \ball (\epsilon_3(H)),
    \end{aligned} $$
    under any disturbance sequence $\{w_k \in \mathcal{W}\}_{t=0}^T$.
\end{lemma}
\begin{proof}
The proof is similar to the proof in \cite{li_online_2021}, except for several notation differences. 
First, it is necessary to verify that $\M(K')\in \mathcal{M}$. 
Indeed, this holds as
$$
\norm{\M^{[i]}(K')}_2 = \norm{(K - K')A_{K'}^{i-1}}_2 \le 2\kappa^3(1-\gamma)^{i-1} = a(1-\gamma)^{i-1},
$$
for all $i \in [H]$.

Then, when applying linear policy $K'$, we have
$$ x_t^{K'} = \sum_{s=1}^t A_{K'}^{s-1} w_{t-s},\quad u_t^{K'} = -\sum_{s=1}^t K' A_{K'}^{s-1} w_{t-s},$$
where $A_{K'} = A - BK'$.

On the other hand, when we implement disturbance-action policy $\pi(K,\M(K'))$, we have
$$ x_t^{\M(K')} = \sum_{s=1}^t \tilde{\Psi}_s^x(\M(K')) w_{t-s}, $$
where 
$$ \tilde{\Psi}_s^x(\M(K')) = A_K^{s-1} + \sum_{j = \max (1,s-H)}^{s-1} A_K^{j-1}B\M^{[s-j]}(K')$$

When $s\le H$, we have
$$
\begin{aligned}
\tilde{\Psi}_s^x(\M(K')) &= A_K^{s-1} + \sum_{j=1}^{s-1} A_K^{j-1}B\M^{[s-j]}(K') \\
&= A_K^{s-1} + \sum_{j=1}^{s-1} A_K^{j-1}B(K-K')(A-BK')^{s-j-1} \\
&= A_K^{s-1} + \sum_{j=1}^{s-1} A_K^{j-1}(A_{K'} - A_K)A_{K'}^{s-j-1} \\
&= A_K^{s-1} + \sum_{j=1}^{s-1} A_K^{j-1}A_{K'}^{s-j} - \sum_{j=1}^{s-1}A_K^jA_{K'}^{s-j-1} \\
& = A_K^{s-1} + A_{K'}^{s-1} - A_K^{s-1} = A_{K'}^{s-1}.
\end{aligned}
$$

When $s>H$,
$$
\begin{aligned}
\tilde{\Psi}_s^x(\M(K')) &= A_K^{s-1} + \sum_{j=s-H}^{s-1} A_K^{j-1}B\M^{[s-j]}(K') \\
&= A_K^{s-1} + \sum_{j=s-H}^{s-1} A_K^{j-1}B(K-K')(A-BK')^{s-j-1} \\
&= A_K^{s-1} + \sum_{j=s-H}^{s-1} A_K^{j-1}A_{K'}^{s-j} - \sum_{j=s-H}^{s-1}A_K^jA_{K'}^{s-j-1} \\
& = A_K^{s-H-1}A_{K'}^H.
\end{aligned}
$$

Hence, 
$$
\begin{aligned}
x_t^{\M(K')} &= \sum_{s=1}^H A_{K'}^{s-1}w_{t-s} + \sum_{s=H+1}^t A_K^{s-H-1}A_{K'}^H w_{t-s} \\
u_t^{\M(K')} &= -Kx_t^{\M(K')} + \sum_{s=1}^H \M^{[s]}(K')w_{t-s} \\
            &= -\sum_{s=1}^H K' A_{K'}^{s-1}w_{t-s} -\sum_{s=H+1}^t K A_{K}^{s-H-1}A_{K'}^H w_{t-s}.
\end{aligned}
$$
Then we can give the bound for $\norm{x_t^{K'} - x_t^{\M(K')}}_{\infty}$ and $\norm{u_t^{K'} - u_t^{\M(K')}}_{\infty}$,

$$
\begin{aligned}
\norm{x_t^{K'} - x_t^{\M(K')}}_{\infty} &\le \norm{x_t^{K'} - x_t^{\M(K')}}_2 = 
\norm{\sum_{s=H+1}^t (A_K^{s-H-1} - A_{K'}^{s-H-1})A_{K'}^Hw_{t-s}}_2\\
& \le \sum_{s=H+1}^t (\norm{A_K^{s-H-1}}_2 + \norm{A_{K'}^{s-H-1}}_2) \norm{A_{K'}^H}_2 \overline{w} \sqrt{n} \\
& \le \sum_{s=H+1}^t 2\kappa^2(1-\gamma)^{s-H-1}\kappa^2(1-\gamma)^H \overline{w}\sqrt{n} \\
& \le 2\kappa^4 \overline{w}\sqrt{n} (1-\gamma)^H/\gamma.
\end{aligned}
$$

$$
\begin{aligned}
\norm{u_t^{K'} - u_t^{\M(K')}}_{\infty} &\le \norm{u_t^{K'} - u_t^{\M(K')}}_2 = 
\norm{\sum_{s=H+1}^t (K A_K^{s-H-1} - K' A_{K'}^{s-H-1})A_{K'}^H w_{t-s}}_2 \\ 
& \le \sum_{s=H+1}^t 2\kappa^3 (1-\gamma)^{s-H-1}\kappa^2(1-\gamma)^H \overline{w}\sqrt{n}\\
& \le 2\kappa^5 \overline{w}\sqrt{n} (1-\gamma)^H/\gamma,
\end{aligned}
$$
under any disturbance sequence $\{ w_k \in \W \}_{k=0}^T$. We can then set $\epsilon_3 (H) := 2\kappa^5 \overline{w}\sqrt{n} (1-\gamma)^H/\gamma = c_3 \sqrt{n}(1-\gamma)^H$ and this finishes the proof.
\end{proof}

Before giving the proof of Theorem \ref{thm: safety}, we need one more corollary that comes directly from Lemma \ref{lemma: epsilon_1 bound} and Lemma \ref{lemma: epsilon_2 bound}.
\begin{corollary} \label{corollary: M(K) in Omega}
Given linear policy $\tilde K \in \K$ that are $\epsilon_0$-strictly safe for $\epsilon_0 \le 0$, we have $\M(\tilde K) \in \Omega^{\epsilon_0 - \epsilon_1 - \epsilon_3}$, where $\M(\tilde K)$ is defined in Lemma \ref{lemma: epsilon_3 bound}.
\end{corollary}
\begin{proof}
    First, note that we use the notation $\tilde{x}_t^{\M(\tilde K)}$ to refer to the surrogate state at time $t$ under the disturbance action policy $\M(\tilde K)$. 
    For the corollary to hold, we need that (i) $\mathring{h}^x(\M(\tilde K)) \cdot \ball^{2 H n} (\overline{w}) \subseteq \X_{\epsilon_0 - \epsilon_1 - \epsilon_3}$ and (ii) $\mathring{h}^u(\M(\tilde K)) \cdot \ball^{2 H n} (\overline{w}) \subseteq \U_{\epsilon_0 - \epsilon_1 - \epsilon_3}$.
    First, we show (i).
    In fact, since $\M(\tilde K)$ is a fixed policy,
    \begin{equation}
    \begin{aligned}
    \label{eqn:hs}
     \mathring{h}^x(\M(\tilde K)) \cdot \ball^{2 H n} (\overline{w}) & = \left\{ \mathring{h}^x(\M(\tilde K)) [w_1^\top \ ...\ w_{2 H}^\top] : w_s \in \mathcal{W}\ \forall s \in [1, 2 H]  \right\},\\
    & = \left\{ \sum_{k=1}^{2H} \mathring{\Psi}_k^x (\M(\tilde K))w_{t-k} : w_s \in \mathcal{W}\ \forall s \in [t-2H, t-1] \right\},\\
    & = \left\{ \tilde{x}_t^{\M(\tilde K)} : w_t \in \mathcal{W}\ \forall t \right\}
    \end{aligned}
    \end{equation}
    and therefore (i) is equivalent to the condition that $\tilde{x}_t^{\M(\tilde K)} \in \X_{\epsilon_* - \epsilon_1 - \epsilon_3}$ under any disturbance sequence $\{w_t \in \mathcal{W}\}$.
    We show this as follows,
    \begin{equation}
    \label{eqn:tils}
    \begin{aligned}
    \tilde{x}_t^{\M(\tilde K)} & = x_t^{\M(\tilde K)} + \left ( \tilde{x}_t^{\M(\tilde K)} - x_t^{\M(\tilde K)} \right)\\
    & \in x_t^{\M(\tilde K)} + \ball^n (\epsilon_1)\\
     &= x_t^{\tilde K} + \left (x_t^{\M(\tilde K)} - x_t^{\tilde K} \right ) + \ball^n (\epsilon_1) \\
    &\subseteq x_t^{\tilde K} +  \ball^n(\epsilon_3) + \ball^n (\epsilon_1)\\
    &\subseteq \X_{\epsilon_0} +  \ball^n(\epsilon_3) + \ball^n (\epsilon_1)\\
    & = \X_{\epsilon_0 - \epsilon_1 - \epsilon_3}.
    \end{aligned}
    \end{equation}
where the first inclusion is due to Lemma \ref{lemma: epsilon_3 bound}, the second inclusion is due to Lemma \ref{lemma: epsilon_1 bound}, the third inclusion is due to Assumption \ref{assumption: existency}, and the fourth inclusion is a property of shrinkage and expansion (see Proposition \ref{proposition: Minkwski subtraction}).
Therefore, combining \eqref{eqn:hs} and \eqref{eqn:tils}, it holds that 
$$
 \mathring{h}^x(\M(\tilde K)) \cdot \ball^{2 H n} (\overline{w}) = \left\{ \tilde{x}_t^{\M(\tilde K)} : w_t \in \mathcal{W}\ \forall t \right\} \subseteq \X_{\epsilon_0 - \epsilon_1 - \epsilon_3}
$$
Condition (ii) holds under the same reasoning.
Since both (i) and (ii) hold, it follows that $\M(\tilde K) \in \Omega^{\epsilon_0 - \epsilon_1 - \epsilon_3}$.
\end{proof}

\noindent \textbf{Proof of Theorem \ref{thm: safety}}
\begin{proof}[Proof of Theorem \ref{thm: safety}] 
    Given $\epsilon, \eta, H$ that satisfy the conditions in Theorem \ref{thm: safety},
    we first prove that $\Omega^{\epsilon}$ is non-empty.
    To show this, we first show that $ \M(K_{ss}) \in \Omega^{\epsilon_* - \epsilon_1 - \epsilon_3}$. This is direct from Corollary \ref{corollary: M(K) in Omega} if we replace $\tilde K$ by $K_{ss}$ and $\epsilon_0$ by $\epsilon_*$. Then by realizing that we choose 
    $\epsilon \le \epsilon_* - \epsilon_1 - \epsilon_3$, we obtain
    $$
        \M(K_{ss}) \in \Omega^{\epsilon_* - \epsilon_1 - \epsilon_3} \subseteq \Omega^{\epsilon},
    $$
which indicates that $\Omega^{\epsilon}$ is non-empty.
Next we show that OGD-BZC is safe. Since it is guaranteed in the algorithm that 
$\M_t \in \Omega^{\epsilon}$ for all $t \ge 0$,  we have 
$$
    \mathring{h}^x(\M_t) \cdot \ball^{2Hn} (\overline{w}) \subseteq \X_{\epsilon}, \quad \forall t \ge 0.
$$

Thus by \eqref{equation: hx, hu} and Lemma \ref{lemma: epsilon_2 bound}, we obtain that
$$\begin{aligned}
    x_{t + 1}^{\A} & = \tilde{x}_{t + 1}^{\A} + (x_{t + 1}^{\A} - \tilde{x}_{t + 1}^{\A})\\
    & \in \tilde{x}_{t + 1}^{\A} + \ball(\epsilon_1)\\
    & = \sum_{k=1}^{2H} \Psi_k^x w_{t+1-k} + \ball(\epsilon_1)\\
    & \subseteq h^x(\M_{t-H+1:t}) \cdot \ball^{2Hn} (\overline{w})  + \ball(\epsilon_1)\\
    & = \left ( \mathring{h}^x(\M_t) + h^x(\M_{t-H+1:t}) - \mathring{h}^x(\M_t) \right ) \cdot \ball^{2Hn} (\overline{w})  + \ball(\epsilon_1) \\
    &\subseteq \mathring{h}^x(\M_t) \cdot \ball^{2Hn} (\overline{w}) + (h^x - \mathring{h}^x) \cdot \ball^{2Hn}(\overline{w})  + \ball(\epsilon_1)\\
    & \subseteq \mathring{h}^x(\M_t) \cdot \ball^{2Hn} (\overline{w}) +\ball^n(\epsilon_2)  + \ball(\epsilon_1)\\
    & \subseteq (\X_{\epsilon})_{-\epsilon_1 -\epsilon_2}\\
    & \subseteq \X_{\epsilon-\epsilon_1 - \epsilon_2},\\
    & \subseteq \X
\end{aligned}
$$
where the first inclusion is due to Lemma \ref{lemma: epsilon_1 bound}, the fourth inclusion is due to Lemma \ref{lemma: epsilon_2 bound}, and the last inclusion is due to the choice $\epsilon \ge \epsilon_1 + \epsilon_2$.
Similarly we can prove that $u_t^{\A}\in \U$. Hence OGD-BZC is safe.
\end{proof}

\newpage

\subsection{Proof related to Theorem \ref{thm: regret bound}}
\label{Appendix: thm2}
\noindent \textbf{Proof of Lemma \ref{lemma: Omega set bound}}
\begin{proof}
We first proof that $ \mathring \Omega(\alpha_1) \subseteq \Omega^{\epsilon}$. For any $\M \in  \mathring \Omega(\alpha_1)$, there exist $\M' \in \Omega$ such that $\M = \alpha_1 \M' + (1-\alpha_1)\M(K_{ss})$
where we use the notation $\M(K_{ss})$ as defined in Lemma \ref{lemma: epsilon_3 bound}. Then,
\begin{equation} 
\label{eqn:hs2}
\begin{aligned}
\mathring{h}^x(\M) \cdot \ball (\overline{w}) &= \mathring{h}^x(\alpha_1 \M' + (1-\alpha_1)\M(K_{ss})) \cdot \ball (\overline{w}) \\
&= \left [ \alpha_1 \mathring{h}^x(\M') + (1-\alpha_1) \mathring{h}^x(\M(K_{ss})) \right ] \cdot \ball (\overline{w}) \\
&\subseteq  \alpha_1 \mathring{h}^x(\M') \cdot \ball (\overline{w}) + (1-\alpha_1) \mathring{h}^x(\M(K_{ss})) \cdot \ball (\overline{w})\\
&\subseteq \alpha_1 \X + (1 -\alpha_1)\X_{\epsilon + r},
\end{aligned}
\end{equation}
where the second equality is due to the fact that $\mathring{h}^x(\M)$ is affine with respect to $\M$. 
In order to show the last inclusion, note that $K_{ss}$ is $\epsilon_*$ strictly safe, and therefore we can apply Corollary \ref{corollary: M(K) in Omega} to get that $\M(K_{ss}) \in \Omega^{\epsilon_* -\epsilon_1 - \epsilon_3} = \Omega^{\epsilon + r}$ and $\mathring{h}^x(\M(K_{ss})) \cdot \ball (\overline{w}) \subseteq \X_{\epsilon + r}$.

Then we claim that $\alpha_1 \X + (1 -\alpha_1)\X_{\epsilon + r} \subseteq \X_{\epsilon}$. This can be verified by following the definition of $\X_{\epsilon}$,

\begin{equation} 
\label{eqn:xs}
\begin{aligned}
    \alpha_1 \X + (1 -\alpha_1)\X_{\epsilon + r} + \ball^n(\epsilon) &= \alpha_1 \X + (1 -\alpha_1)(\X_{\epsilon + r} + \frac{1}{1 - \alpha_1}\ball^n(\epsilon))\\
    &= (1 - \frac{\epsilon}{\epsilon+ r})\X + \frac{\epsilon}{\epsilon+r} (\X_{\epsilon + r} + \ball^n(\epsilon+ r)) \\
    &\subseteq (1 - \frac{\epsilon}{\epsilon+ r})\X + \frac{\epsilon}{\epsilon+ r}\X\\
    & \subseteq \X.
\end{aligned}
\end{equation}
where the last inclusion is due to the convexity of $\X$. 
Combining \eqref{eqn:hs2} and \eqref{eqn:xs}, we have that $\mathring{h}^x(\M) \cdot \ball (\overline{w}) \subseteq \X_{\epsilon}$. Similarly, we can also derive $\mathring{h}^u(\M) \cdot \ball (\overline{w}) \subseteq \U_{\epsilon}$. Hence we can conclude that $ \mathring \Omega(\alpha_1) \subseteq \Omega^{\epsilon}$.

Next we prove $ \mathring \Omega(\alpha_2) \supseteq \Omega^{-\epsilon_1 - \epsilon_3}$. This is equivalent to $\Omega^{-\epsilon_1 - \epsilon_3} \subseteq \alpha_2\Omega + (1-\alpha_2)\M(K_{ss})$. 
Rearranging this inclusion yields
$$
\begin{aligned}
\Omega^{-\epsilon_1 - \epsilon_3} \subseteq \alpha_2\Omega + (1-\alpha_2)\M(K_{ss}) \quad  & \Longleftrightarrow \quad \Omega^{-\epsilon_1 - \epsilon_3} + (\alpha_2 - 1)\M(K_{ss}) \subseteq \alpha_2\Omega \\
& \Longleftrightarrow \quad \frac{1}{\alpha_2} \Omega^{-\epsilon_1 -\epsilon_3} + (1 - \frac{1}{\alpha_2})\M(K_{ss}) \subseteq \Omega
\end{aligned}
$$
We will show this directly.
First note that, for any $\M \in \frac{1}{\alpha_2} \Omega^{-\epsilon_1 -\epsilon_3} + (1 - \frac{1}{\alpha_2})\M(K_{ss})$, there exist $\M' \in \Omega^{-\epsilon_1-\epsilon_3}$ such that $\M = \frac{1}{\alpha_2} \M' + (1 - \frac{1}{\alpha_2})\M(K_{ss})$.
Therefore, we will show that $\M \in \Omega$, by showing that both (i) $\mathring{h}^x(\M) \cdot \ball (\overline{w}) \in \X$ and (ii)  $\mathring{h}^u(\M) \cdot \ball (\overline{w}) \in \U$.
We show (i) as follows
$$
\begin{aligned}
    \mathring{h}^x(\M) \cdot \ball (\overline{w}) &= \mathring{h}^x(\frac{1}{\alpha_2} \M' + (1 - \frac{1}{\alpha_2})\M(K_{ss})) \cdot \ball (\overline{w}) \\
    &= \left [ \frac{1}{\alpha_2} \mathring{h}^x(\M') + (1-\frac{1}{\alpha_2})\mathring{h}^x(\M(K_{ss})) \right ] \cdot \ball(\overline{w})  \\
    &\subseteq \frac{1}{\alpha_2} \mathring{h}^x(\M') \cdot \ball(\overline{w}) + (1 -\frac{1}{\alpha_2}) \mathring{h}^x(\M(K_{ss})) \cdot \ball(\overline{w}) \\
    &\subseteq \frac{1}{\alpha_2} (\X + \ball^n(\epsilon_1 + \epsilon_3)) + (1 - \frac{1}{\alpha_2}) \X_{\epsilon + r} \\
    & = \frac{1}{\alpha_2}\X + (1 - \frac{1}{\alpha_2}) (\X_{\epsilon + r} + \frac{1}{\alpha_2 - 1}\ball^n(\epsilon_1 + \epsilon_3))\\
    & = \frac{1}{\alpha_2}\X + (1 - \frac{1}{\alpha_2}) (\X_{\epsilon + r} + \ball^n(\epsilon + r))\\
    & \subseteq \frac{1}{\alpha_2}\X + (1 - \frac{1}{\alpha_2})\X \subseteq \X, \\
\end{aligned}
$$
where the second inclusion follows the same reasoning as in \eqref{eqn:hs2}.
Point (ii) can be shown similarily. This indicates that $\M \in \Omega$. Therefore, we have $ \mathring \Omega(\alpha_2) \supseteq \Omega^{-\epsilon_1 - \epsilon_3}$.
\end{proof}

\noindent \textbf{Proof of Lemma \ref{lemma: M_ap}}
\begin{proof}
The statement that $\M_{\text{ap}} \in \Omega^{-\epsilon_1 - \epsilon_3}$ can directly derived from Corollary \ref{corollary: M(K) in Omega} by noticing that $K^* \in \K$ and thus $0$-strictly safe. Then we will bound $J_T (\M_{\text{ap}}) - J_T(K^*)$. Let $x_t,u_t$ and $x_t^*, u_t^*$ denotes the states/response by using policy $\pi(K,\{\M_{\text{ap}}\})$ and linear policy $K^*$, then
$$
\begin{aligned}
  J_T (\M_{\text{ap}}) - J_T(K^*) &= \sum_{t=0}^T [c_t(x_t,u_t) - c_t (x_t^*, u_t^*)]\\
  &\le \sum_{t=0}^TGb(\norm{x_t - x_t^*}_2 + \norm{u_t - u_t^*}_2)\\
  &\le \sum_{t=0}^TGb(\sqrt{n}\norm{x_t - x_t^*}_{\infty} + \sqrt{m}\norm{u_t - u_t^*}_{\infty}) \\
  &\le \sum_{t=0}^TGb \sqrt{mn} (\norm{x_t - x_t^*}_{\infty} + \norm{u_t - u_t^*}_{\infty}) \\
  &\le 2TGb \sqrt{mn} \cdot \epsilon_3,
\end{aligned}
$$
where the last equality comes from Lemma \ref{lemma: epsilon_3 bound}.
Then by Lemma \ref{lemma: x_t} and Assumption \ref{assumption: loss function}, it is direct to see that 
$$
\begin{aligned}
    \sum_{t=0}^T \mathring{f}_t(\M_{\text{ap}})-\min_{K \in \K} J_T(K)
    &= \sum_{t=0}^T \mathring{f}_t(\M_{\text{ap}}) - J_T(\M_{\text{ap}})+ J_T (\M_{\text{ap}}) - J_T(K^*)\\
    &\le \sum_{t=0}^T Gb (\norm{x_t - \tilde{x}_t}_2 + \norm{u_t - \tilde{u}_t}_2)  + 2TGb \sqrt{mn} \cdot \epsilon_3 \\
    &\le \sum_{t=0}^T 2Gb \sqrt{mn} \cdot \epsilon_1 + 2TGb \sqrt{mn} \cdot \epsilon_3\\
    &\le 2 T Gb \sqrt{mn} \cdot (\epsilon_1 + \epsilon_3) \\
    &\le \mathcal{O}(TmnH(\epsilon_1 + \epsilon_3)),
\end{aligned}
$$
when $H \ge \frac{\log (2\kappa^2)}{\log((1-\gamma)^{-1})}$.
\end{proof}

\noindent \textbf{Proof of Lemma \ref{lemma: Performance gap}}
\begin{proof}
We will prove Lemma \ref{lemma: Performance gap} by giving three separate bounds as follows:
$$
\begin{aligned}
     J_T(\A)- \min _{\M \in \Omega^{\epsilon}} \sum_{t=0}^T \mathring{f}_t(\M) &= \\
     \underbrace{J_T (\A) - \sum_{t=0}^T f_t(\M_{t-H:t})}_{\text{Part i}} &+ 
     \underbrace{\sum_{t=0}^T f_t(\M_{t-H:t}) - \sum_{t=0}^T \mathring{f}_t\left(\M_t\right)}_{\text{Part ii}} + 
     \underbrace{\sum_{t=0}^T \mathring{f}_t\left(\M_t\right) -\min _{\M \in \Omega^{\epsilon}} \sum_{t=0}^T \mathring{f}_t(\M)}_{\text{Part iii}}.
\end{aligned}
$$

For Part ii, by Lemma \ref{lemma: bound on f_t}, we have
$$ |f_t(\M_{t-H:t}) - \mathring{f}_t(\M_t)| \le 2Gb \left( \sqrt n \bar w (1 + \kappa) \kappa^2\kappa_B \sum_{i=1}^H(1-\gamma)^{i-1} \sum_{j=1}^H \norm{\M_{t-i}^{[j]}- \M_t^{[j]}}_2 \right ).$$

Since $\M_t$ is generated by the OGD-BZC algorithm with fixed step size $\eta$, by Lemma \ref{lemma: bound on f_t} we have

$$\norm{\M_t- \M_{t-i}}_F \le \sum_{s=1}^{i} \norm{\M_{t-s+1}- \M_{t-s}}_F \le 
\sum_{s=1}^{i} \eta \norm{\nabla \mathring{f}_{t-s}(\M_{t-s})}_F \le i\eta G_f,$$
and $Gb \le \mathcal{O}(\sqrt{mn}H), Gf \le \mathcal{O}(\sqrt{mn^2H^3})$. By plugging in the above result, we have for all $t>0$,
$$ 
\begin{aligned}
|f_t(\M_{t-H:t}) - \mathring{f}_t(\M_t)| &\le 2Gb \left( \sqrt n \bar w (1 + \kappa) \kappa^2\kappa_B \sqrt{H} \sum_{i=1}^H(1-\gamma)^{i-1} \norm{\M_t- \M_{t-i}}_F  \right ) \\
& \le 2Gb \left( \sqrt n \bar w (1 + \kappa) \kappa^2\kappa_B \sqrt{H} \sum_{i=1}^H(1-\gamma)^{i-1} i \eta G_f  \right ) \\
& \le 2Gb  \sqrt{n} \bar w (1 + \kappa) \kappa^2\kappa_B \sqrt{H} \eta G_f \frac{1}{\gamma^2} \\
& \le \mathcal{O} ({\sqrt{mn}H \cdot \sqrt{mn^2H^3} \cdot \eta \sqrt{nH}}) \\
& = \mathcal{O} (mn^2H^3\eta).
\end{aligned}
$$

Thus $\text{Part ii} \le \sum_{t=0}^T |f_t(\M_{t-H:t}) - \mathring{f}_t(\M_t)| \le \mathcal{O}(Tmn^2H^3\eta) $. 

For Part i, let $x_t,u_t$ denotes the states/response by running the algorithm, then by Assumption \ref{assumption: loss function} we have

$$
\begin{aligned}
  J_T (\A) - \sum_{t=0}^T f_t(\M_{t-H:t}) &= \sum_{t=0}^T [c_t(x_t,u_t) - c_t (\tilde{x}_t, \tilde{u}_t)]\\
  &\le \sum_{t=0}^TGb(\norm{x_t - \tilde{x}_t}_2 + \norm{u_t - \tilde{u}_t}_2)\\
  &\le \sum_{t=0}^TGb(\sqrt{n}\norm{x_t - \tilde{x}_t}_{\infty} + \sqrt{m}\norm{u_t - \tilde{u}_t}_{\infty}) \\
  &\le \sum_{t=0}^TGb \sqrt{mn} (\norm{x_t - \tilde{x}_t}_{\infty} + \norm{u_t - \tilde{u}_t}_{\infty}) \\
  &\le 2TGb \sqrt{mn} \cdot \epsilon_1 \le \mathcal{O} (T \cdot \sqrt{mn} H \cdot \sqrt{mn} H (1-\gamma)^H) \\
  & = \mathcal{O} (T mn H^2 (1-\gamma)^H).
\end{aligned}
$$

For Part iii analysis is a standard gradient descent analysis that can be found in \cite{hazan_introduction_2023}. From the literature, we have
$$
\sum_{t=0}^T \mathring{f}_t\left(\M_t\right) -\min _{\M \in \Omega^{\epsilon}} \sum_{t=0}^T \mathring{f}_t(\M) \le \frac{\delta^2}{\eta} + G_f^2 T \eta.
$$

Also by Lemma \ref{lemma: bound on f_t}, we have $G_f := \sup_{t \in [T],\M \in \Omega^{\epsilon}} ||\nabla \mathring{f}_t(\M)||_F \le \mathcal{O} (\sqrt{mn^2H^3})$. It only remains to show that $\delta := \sup_{\M, \tilde{\M} \in \Omega^{\epsilon}} ||\M - \tilde{\M}||_F \le 2 \sqrt{m}a/\gamma$. In fact,
$$
\begin{aligned}
    \delta &= \sup_{\M, \tilde{\M} \in \Omega^{\epsilon}} \norm{\M - \tilde{\M}}_F \le \sup_{\M, \tilde{\M} \in \mathcal{M}} \norm{\M - \tilde{\M}}_F \\
    &= \sup_{\M, \tilde{\M} \in \mathcal{M}} \sum_{i=1}^{H}\norm{\M^{[i]} - \tilde{\M}^{[i]}}_F \\
    &\le \sup_{\M, \tilde{\M} \in \mathcal{M}} \sqrt{m} \sum_{i=1}^{H}\norm{\M^{[i]} - \tilde{\M}^{[i]}}_2 \\
    &\le \sup_{\M, \tilde{\M} \in \mathcal{M}} \sqrt{m} \sum_{i=1}^{H}\norm{\M^{[i]}}_2 + \norm{\tilde{\M}^{[i]}}_2 \\
    &\le 2\sqrt{m} \sum_{i=1}^H a(1-\gamma)^i \le 2\sqrt{m}a/\gamma.
\end{aligned}
$$

Combine the bound for Part i, Part ii, and Part iii, we will have 
$$
\begin{aligned}
J_T(\A)- \min _{\M \in \Omega^{\epsilon}} \sum_{t=0}^T \mathring{f}_t(\M) &\le \mathcal{O}( T mn H^2 (1-\gamma)^H + Tmn^2H^3\eta + Tmn^2H^3\eta + \frac{m}{\eta}) \\
&= \mathcal{O}(T mn H^2 (1-\gamma)^H + Tmn^2H^3\eta + \frac{m}{\eta}).
\end{aligned}
$$

This finishs the proof.
\end{proof}

\newpage

\subsection{Helping Lemmas}
\label{Appendix: helping}
The lemmas in this section are fundamental analysis steps and are closely drawn from \cite{li_online_2021} and \cite{agarwal_online_2019}, with slight notation differences. We include them here as auxiliary lemmas, which will be referenced when proving other lemmas.

\begin{lemma}[Bound on $\Psi_k^x$ and $\Psi_k^u$] \label{lemma: Psi_k} \label{Bound on disturbance-to-state/response matrix}
Suppose $\M_t \in \mathcal{M}$ for all $t$, we have
$$
\begin{aligned}
    &\norm{\Psi_k^x(\M_{t-H:t-1})}_2 \le \kappa^2 (1-\gamma)^{k-1} \one_{k \le H} + a\kappa^2 \kappa_B H (1-\gamma)^{k-2} \one_{(k\ge 2)} \\
    &\norm{\Psi_k^u(\M_{t-H:t})}_2 \le (a + \kappa^3)(1-\gamma)^{k-1}\one_{(k\le H)} + a \kappa^3\kappa_B H (1-\gamma)^{k-2} \one_{(k\ge 2)}.
\end{aligned}
$$
\begin{proof}
    By Definition \ref{def: DAC} and Proposition \ref{prop: DAC}, we have
    $$
    \begin{aligned}
        \norm{\Psi_k^x(\M_{t-H:t-1})}_2 &= \norm{A_k^{k-1} \one_{(k\le H)} + \sum_{i=1}^H A_K^{i-1} B M_{t-i}^{[k-i]} \one_{(1 \le k-i \le H )} }_2 \\
        &\le \norm{A_k^{k-1}}_2 \one_{(k \le H)} + \sum_{i=1}^H \norm{A_K^{i-1}}_2 \norm{B}_2 \norm{\M_{t-i}^{[k-i]}}_2 \one_{(1 \le k-i \le H )} \\
        &\le \kappa^2 (1-\gamma)^{k-1} \one_{k \le H} + \sum_{i=1}^H \kappa^2 (1-\gamma)^{i-1} \kappa_B \norm{\M_{t-i}^{[k-i]}}_2 \one_{(1\le k-i \le H)} \\
        &\le \kappa^2 (1-\gamma)^{k-1} \one_{k \le H} + \sum_{i=1}^H \kappa^2 (1-\gamma)^{i-1} \kappa_B \cdot a (1-\gamma)^{k-i-1} \one_{(1\le k-i \le H)} \\
        &\le \kappa^2 (1-\gamma)^{k-1} \one_{k \le H} + \sum_{i=1}^H a\kappa^2 \kappa_B (1-\gamma)^{k-2} \one_{(1\le k-i \le H)} \\
        &\le \kappa^2 (1-\gamma)^{k-1} \one_{k \le H} + a\kappa^2 \kappa_B H (1-\gamma)^{k-2} \one_{(k\ge 2)}. 
    \end{aligned}
    $$
    Further for $\Psi_k^u$ we have
    $$
    \begin{aligned}
        \norm{\Psi_k^u(\M_{t-H:t})}_2 &= \norm{\M_t^{[k]} \one_{(k \le H)} - K \Psi_k^x(\M_{t-H:t-1})}_2 \\
        &\le \norm{\M_t^{[k]}}_2 \one_{(k \le H)} + \norm{K}_2 \norm{\Psi_k^x(\M_{t-H:t-1})}_2 \\
        &\le a(1-\gamma)^{k-1} \one_{(k\le H)} + \kappa \left(\kappa^2 (1-\gamma)^{(k-1)} \one_{k \le H} + a\kappa^2 \kappa_B H (1-\gamma)^{k-2} \one_{(k\ge 2)} \right) \\
        &\le (a + \kappa^3)(1-\gamma)^{k-1}\one_{(k\le H)} + a \kappa^3\kappa_B H (1-\gamma)^{k-2} \one_{(k\ge 2)}.
    \end{aligned}
    $$
\end{proof}

\end{lemma}

\begin{lemma}[Bound on $\norm{x_t}_2$ and $\norm{u_t}_2$]\label{lemma: x_t} 
Suppose $\M_t \in \mathcal{M}$ for all $t$, we have
$$ \max(\norm{x_t}_2, \norm{\tilde{x}_t}_2) \le b_x, \quad \max(\norm{u_t}_2, \norm{\tilde{u}_t}_2) \le b_u,  $$
where $$b_x := \frac{\overline{w}\sqrt{n} (\kappa^2 + a\kappa^2 \kappa_B H)}{(1-\kappa^2 (1-\gamma)^H)\gamma}, \quad b_u := \kappa b_x+ \overline{w}\sqrt{n} \cdot a /\gamma. $$ 
If we further define $b := \max(b_x, b_u)$ then we have 
$$\max(\norm{x_t}_2, \norm{\tilde{x}_t}_2,\norm{u_t}_2, \norm{\tilde{u}_t}_2) \le b$$
When $H \ge \frac{\log (2\kappa^2)}{\log((1-\gamma)^{-1})}$, we have 
$$b \le 2 \overline{w}\sqrt{n} (\kappa^3 + a\kappa^3 \kappa_B H)/\gamma + \overline{w}\sqrt{n}  \cdot a /\gamma \le \mathcal{O}(\sqrt{mn}H) $$

\begin{proof}
We first bound $\norm{\tilde{x}_t}_2$ and $\norm{\tilde{u}_t}_2$. By Lemma \ref{Bound on disturbance-to-state/response matrix} we have,
$$
\begin{aligned}
    \norm{\tilde{x}_t}_2 &= \norm{\sum_{k=1}^{2H}\Psi_k^x w_{t-k}}_2 \le  \overline{w} \sqrt{n} \sum_{k=1}^{2H} \norm{\Psi_k^x}_2 \\
    &\le \overline{w} \sqrt{n} \sum_{k=1}^{2H} (\kappa^2 (1-\gamma)^{k-1} \one_{k \le H} + a\kappa^2 \kappa_B H (1-\gamma)^{k-2} \one_{(k\ge 2)})   \\
    &\le \overline{w} \sqrt{n} (\kappa^2 + a\kappa^2 \kappa_B H)/\gamma \le b_x,
\end{aligned}
$$ 
and 
$$
\begin{aligned}
    \norm{\tilde{u}_t}_2 &= \norm{\sum_{k=1}^{2H} \Psi_k^u w_{t-k}}_2 \le \overline{w} \sqrt{n} \sum_{k=1}^{2H} \norm{\Psi_k^u}_2 \\
    &\le \overline{w} \sqrt{n} \sum_{k=1}^{2H} \left ( (a + \kappa^3)(1-\gamma)^{k-1}\one_{(k\le H)} + a \kappa^3\kappa_B H (1-\gamma)^{k-2} \one_{(k\ge 2)} \right ) \\
    &\le \overline{w} \sqrt{n} \left ( a + \kappa^3 + a\kappa^3 \kappa_B H \right )/\gamma \le \kappa b_x+ \overline{w}\sqrt{n} \cdot a /\gamma = b_u.
\end{aligned}
$$ 
Further, based on the fact that $\norm{x_0} = \norm{u_0} = 0$ and by induction, we can prove that
$$ \norm{x_t}_2 = \norm{A_K^H x_{t-H} + \tilde{x}_t}_2 \le \kappa^2 (1-\gamma)^H b_x + \overline{w} \sqrt{n} (\kappa^2 + a\kappa^2 \kappa_B H)/\gamma \le b_x,
$$
and
$$
\begin{aligned}
\norm{u_t}_2 & =\norm{-K x_t+\sum_{i=1}^H \M_t^{[i]} w_{t-i}}_2 \le \norm{-K x_t}_2 +\sum_{i=1}^H \norm{\M_t^{[i]} w_{t-i}}_2 \\ 
& \le \kappa b_x+\sum_{i=1}^H \sqrt{n} \cdot a (1-\gamma)^{i-1} \overline{w} \\ 
& \le \kappa b_x+ \overline{w}\sqrt{n} \cdot a /\gamma = b_u.
\end{aligned}
$$

\end{proof}
\end{lemma}

\begin{lemma}[Bound on $f_t$ and $\nabla \mathring{f}_t$] \label{lemma: bound on f_t}
For any $\M_t, \tilde{\M_t} \in \mathcal{M}$ for all $t$, we have
$$ |f_t(\M_{t-H:t-1}) - f_t(\tilde{\M}_{t-H:t-1})| \le 2Gb \left( \sqrt n \bar w (1 + \kappa) \kappa^2\kappa_B \sum_{i=1}^H(1-\gamma)^{i-1} \sum_{j=1}^H \norm{\M_{t-i}^{[j]}-\tilde \M_{t-i}^{[j]}}_2 \right ).$$

Further, for any $t>0$ and $\M_t \in \mathcal{M}$, we have $\norm{\nabla \mathring{f}_t (\M_t)}_F \le G_f := 2G\tilde{b} \sqrt n \bar w (1 + \kappa) \kappa^2\kappa_B \sqrt{H}/\gamma$.
When $H \ge \frac{\log(2\kappa^2)}{\log((1-\gamma)^{-1})}$, we have $\tilde{b} \le 2 \overline{w}\sqrt{n} (\kappa^3 + \tilde{a}\kappa^3 \kappa_B H)/\gamma + \overline{w}\sqrt{n} \cdot \tilde{a} /\gamma$, and the result simpifies to $G_f \le \mathcal{O}(\sqrt{m n^2 H^3})$.

\begin{proof}
    Let $\tilde{x}$ and $\tilde{\tilde{x}}$ be the surrogate state generated by $\pi (K, \{\M_t\})$ and $\pi (K, \{\tilde{\M}_t\})$, we have
    $$
    \begin{aligned}
	\norm{\tilde x_t- \tilde{\tilde{x}}_t}_2& =\norm{\sum_{k=1}^{2H} (\Psi_{k}^x( \M_{t-H:t-1})-\Psi_{k}^x(\tilde{\M}_{t-H:t-1}))w_{t-k}}_2 \\
	& \le \sum_{k=1}^{2H} \norm{\Psi_{k}^x(\M_{t-H:t-1})-\Psi_{k}^x(\tilde{\M}_{t-H:t-1})}_2\sqrt n \bar w\\
	& \le \sqrt n \bar w\sum_{k=1}^{2H} \norm{\sum_{i=1}^H A_K^{i-1}B (\M_{t-i}^{[k-i]}-\tilde \M_{t-i}^{[k-i]})\one_{(1\le k-i \le H)}}_2\\
	&\le \sqrt n \bar w\sum_{k=1}^{2H}\sum_{i=1}^H \kappa^2 (1-\gamma)^{i-1}\kappa_B \norm{\M_{t-i}^{[k-i]}-\tilde \M_{t-i}^{[k-i]}}_2\one_{(1\le k-i \le H)}\\
	& = \sqrt n \bar w\kappa^2\kappa_B \sum_{i=1}^H(1-\gamma)^{i-1} \sum_{j=1}^H \norm{\M_{t-i}^{[j]}-\tilde \M_{t-i}^{[j]}}_2,
    \end{aligned}
    $$
and let $\tilde u_t, \tilde{\tilde{u}}$ be the corresponding surrogate response, we have
    $$
    \begin{aligned}
	\norm{\tilde u_t-\tilde{\tilde{u}}_t}_2&=\norm{-K \tilde x_t + K \tilde{\tilde{x}}_t+ \sum_{i=1}^H \M_t^{[i]}w_{t-i}- \sum_{i=1}^H \tilde \M_t^{[i]}w_{t-i}}_2 \\
	&\le \kappa 	\norm{\tilde x_t-\tilde{\tilde{x}}_t}_2 + \sqrt n  \bar w \sum_{i=1}^H \norm{\M_t^{[i]}-\tilde \M_t^{[i]}}_2\\
        & = \sqrt n \bar w\kappa^3 \kappa_B \sum_{i=1}^H(1-\gamma)^{i-1} \sum_{j=1}^H \norm{\M_{t-i}^{[j]}-\tilde \M_{t-i}^{[j]}}_2 + \sqrt n  \bar w \sum_{i=1}^H \norm{\M_t^{[i]}-\tilde \M_t^{[i]}}_2.\\
    \end{aligned}
    $$
Then combine Assumption \ref{assumption: loss function} and Lemma \ref{lemma: x_t}, we get 
    \begin{equation} \label{eq: f_t(M1) - f_t(M2)}
    \begin{aligned}
	|f_t(\M_{t-H:t-1}) - f_t(\tilde{\M}_{t-H:t-1})| &= |c_t(\tilde{x}_t, \tilde{u}_t) - c_t(\tilde{\tilde{x}}_t, \tilde{\tilde{u}}_t) | \\
        & \le \norm{\nabla c_t}_2 \left (\norm{\tilde{x}_t - \tilde{\tilde{x}}_t}_2 + \norm{\tilde{u}_t - \tilde{\tilde{u}}_t}_2 \right ) \\
        & \le Gb \left (\norm{\tilde{x}_t - \tilde{\tilde{x}}_t}_2 + \norm{\tilde{u}_t - \tilde{\tilde{u}}_t}_2 \right ) \\
        & \le Gb \left( \sqrt n \bar w (1 + \kappa) \kappa^2\kappa_B \sum_{i=1}^H(1-\gamma)^{i-1} \sum_{j=1}^H \norm{\M_{t-i}^{[j]}-\tilde \M_{t-i}^{[j]}}_2 +  \sqrt n  \bar w \sum_{i=1}^H \norm{\M_t^{[i]}-\tilde \M_t^{[i]}}_2 \right) \\
        & \le 2Gb \left( \sqrt n \bar w (1 + \kappa) \kappa^2\kappa_B \sum_{i=1}^H(1-\gamma)^{i-1} \sum_{j=1}^H \norm{\M_{t-i}^{[j]}-\tilde \M_{t-i}^{[j]}}_2 \right ).
    \end{aligned}
    \end{equation}
where the last inequality uses the fact that $\kappa, \kappa_B > 1$. 

Further, we can establish a bound for $\norm{\nabla \mathring{f}_t (\M)}_F$ when $\M \in \mathcal{M}$. To achieve this, we define a slightly larger set $\tilde{\mathcal{M}} := \{ \M: \norm{\M^{[i]}}_2 \le \tilde{a}(1-\gamma)^{i-1}, \forall i \in [H] \}$, where $\tilde{a} = 2a$, for some small $\delta >0$. Then, by following the same analysis as in Lemma \ref{lemma: Psi_k} and Lemma \ref{lemma: x_t}, we will obtain a result similar to \eqref{eq: f_t(M1) - f_t(M2)}. Specifically, when $\M \in \mathcal{M}$ and $\M + \Delta \M \in \tilde{\mathcal{M}}$, we have
$$
| \mathring{f}_t(\M + \Delta \M) - \mathring{f}(\M)| \le 2G \tilde{b} \left( \sqrt{n} \bar{w} (1 + \kappa) \kappa^2\kappa_B \sum_{i=1}^H(1-\gamma)^{i-1} \sum_{j=1}^H \norm{\Delta \M^{[j]}}_2 \right ),
$$
where $\tilde{b} \le 2 \overline{w}\sqrt{n} (\kappa^3 + \tilde{a}\kappa^3 \kappa_B H)/\gamma + \overline{w}\sqrt{mn} \cdot \tilde{a} /\gamma \le \mathcal{O}(\sqrt{n}H)$, when $H \ge \frac{\log(2\kappa^2)}{\log((1-\gamma)^{-1})}$.

Let $\M \in \mathcal{M}$ and we have

$$
\begin{aligned}
    \norm{\nabla \mathring{f}_t (\M)}_F &=  \sup_{\Delta \M \neq 0; \M, \M + \Delta \M \in \tilde{\mathcal{M}}} \frac{  \langle \nabla \mathring{f}_t (\M), \Delta \M \rangle}{\norm{\Delta \M}_F} \\
    & \le  \sup_{\Delta \M \neq 0; \M, \M + \Delta \M \in \tilde{\mathcal{M}}} \frac{\mathring{f}_t (\M + \Delta \M) - \mathring{f}_t (\M)}{\norm{\Delta \M}_F} \\
    & \le \sup_{\Delta \M \neq 0; \M, \M + \Delta \M \in \tilde{\mathcal{M}}} \frac{2G\tilde{b} \left( \sqrt n \bar w (1 + \kappa) \kappa^2\kappa_B \sum_{i=1}^H(1-\gamma)^{i-1} \sqrt{H} \norm{\Delta \M}_F \right )}{\norm{\Delta \M}_F} \\
    & \le 2G\tilde{b}  \sqrt n \bar w (1 + \kappa) \kappa^2\kappa_B \sqrt{H}/\gamma
\end{aligned}
$$

Then we have $\norm{\nabla \mathring{f}_t (\M)}_F \le 2G\tilde{b} \left( \sqrt n \bar w (1 + \kappa) \kappa^2\kappa_B \sqrt{H}/\gamma \right)$.



Setting $G_f = 2G\tilde{b}  \sqrt n \bar w (1 + \kappa) \kappa^2\kappa_B \sqrt{H}/\gamma $, we have $\norm{\nabla \mathring{f}_t (\M)}_F \le G_f $. This finishes the proof.
     
\end{proof}
\end{lemma}

\begin{lemma} \label{lemma: Psi difference}
    Consider any $\M_t, \tilde{\M_t} \in \mathcal{M}$ for all $t$, we have

    $$
    \begin{aligned}
    &\sum_{k=1}^{2H} \norm{\Psi_k^x(\M_{t-H:t-1}) - \Psi_k^x(\tilde{\M}_{t-H:t-1})}_2 \le  \kappa^2 \kappa_B \sqrt{H} \sum_{i=1}^{H}  (1 - \gamma)^{i-1} \norm{ \M_{t-i} - \tilde{\M}_{t-i} }_F,\\
    &\sum_{k=1}^{2H} \norm{\Psi_k^u(\M_{t-H:t}) - \Psi_k^u(\tilde{\M}_{t-H:t})}_2 
    \le \sqrt{H} \norm{\M_t -\tilde{\M}_t}_F + \kappa^3 \kappa_B \sqrt{H} \sum_{i=1}^H (1-\gamma)^{i-1} \norm{ \M_{t-i} - \tilde{\M}_{t-i} }_F.
     \end{aligned}
    $$
\begin{proof}
    We first provide a bound for $\norm{\Psi_k^x(M_{t-H:t-1}) - \Psi_k^x(\tilde{M}_{t-H:t-1})}_2$,
    $$ 
    \begin{aligned}
    \norm{ \Psi_k^x (\M_{t+k:t-1}) - \Psi_k^x (\tilde{\M}_{t+k:t-1}) }_2 &= \norm{ \sum_{i=1}^{H} A^{i-1}_{k} B \M^{[k-i]}_{t-i} \one_{(1 \le k-i \le H)} - \sum_{i=1}^{H} A^{i-1}_{k} B \tilde{\M}_{t-i}^{[k-i]} \one_{(1 \le k-i \le H)} }_2 \\
    &= \norm{ \sum_{i=1}^{H} A^{i-1}_{k} B \left( \M_{t-i}^{[k-i]} - \tilde{\M}_{t-i}^{[k-i]} \right) \one_{(1 \le k-i \le H)} }_2 \\
    &\le \sum_{i=1}^{H} \norm{ A^{i-1}_{k} }_2 \norm{ B }_2 \norm{ \M^{[k-i]}_{t-i} - \tilde{\M}^{[k-i]}_{t-i} }_2 \one_{(1 \le k-i \le H)} \\
    &\le \sum_{i=1}^{H} \kappa^2 (1 - \gamma)^{i-1} \kappa_B \norm{ \M^{[k-i]}_{t-i} - \tilde{\M}^{[k-i]}_{t-i} }_2 \one_{(1 \le k-i \le H)}.
    \end{aligned}
    $$
    Then 
    $$ 
    \begin{aligned}
    \sum_{k=1}^{2H} \norm{\Psi_k^x(M_{t-H:t-1}) - \Psi_k^x(\tilde{M}_{t-H:t-1})}_2 
    &\le \sum_{k=1}^{2H} \sum_{i=1}^{H} \kappa^2 (1 - \gamma)^{i-1} \kappa_B \norm{ \M^{[k-i]}_{t-i} - \tilde{\M}^{[k-i]}_{t-i} }_2 \one_{(1 \le k-i \le H)} \\
    &=\kappa^2 \kappa_B \sum_{i=1}^{H}  (1 - \gamma)^{i-1} \sum_{k=1}^{2H}\norm{ \M^{[k-i]}_{t-i} - \tilde{\M}^{[k-i]}_{t-i} }_2 \one_{(1 \le k-i \le H)} \\
    & \le \kappa^2 \kappa_B \sqrt{H} \sum_{i=1}^{H}  (1 - \gamma)^{i-1} \norm{ \M_{t-i} - \tilde{\M}_{t-i} }_F.
    \end{aligned}
    $$
    
    Next we provide a bound for $\norm{\Psi_k^u(\M_{t-H:t-1}) - \Psi_k^u(\tilde{\M}_{t-H:t-1})}_2$,
    $$ 
    \begin{aligned}
    \norm{ \Psi_k^u (\M_{t+k:t-1}) - \Psi_k^u (\tilde{\M}_{t+k:t-1}) }_2 
    &\le \norm{\M_t^{[k]}- \tilde{\M}_t^{[k]}}_2 \one_{(k\le H)} + \norm{ \sum_{i=1}^{H}K A^{i-1}_{k} B \left( \M_{t-i}^{[k-i]} - \tilde{\M}_{t-i}^{[k-i]} \right) }_2 \one_{(1 \le k-i \le H)} \\
    &\le \norm{\M_t^{[k]}- \tilde{\M}_t^{[k]}}_2 \one_{(k\le H)} + \sum_{i=1}^{H} \kappa^3 (1 - \gamma)^{i-1} \kappa_B \norm{ \M^{[k-i]}_{t-i} - \tilde{\M}^{[k-i]}_{t-i} }_2 \one_{(1 \le k-i \le H)}
    \end{aligned}
    $$
    Therefore,
    $$ 
    \begin{aligned}
    \sum_{k=1}^{2H} \norm{\Psi_k^u(\M_{t-H:t}) - \Psi_k^u(\tilde{\M}_{t-H:t})}_2 
    &\le \sqrt{H} \norm{\M_t -\tilde{\M}_t}_F + \kappa^3 \kappa_B \sqrt{H} \sum_{i=1}^H (1-\gamma)^{i-1} \norm{ \M_{t-i} - \tilde{\M}_{t-i} }_F. \\
    \end{aligned}
    $$

\end{proof}
\end{lemma}

\end{document}